\documentclass[a4paper]{article}

\usepackage{latexsym}

\newtheorem{thm}{Theorem}[section]

\newtheorem{lem}[thm]{Lemma}
\newtheorem{cor}[thm]{Corollary}

\def\qed{{$\Box$}\medskip}

\def\implies{\Rightarrow}

\def\a{\alpha}  \def\g{\gamma} \def\d{\delta}
\def\e{\varepsilon} \def\f{\varphi} \def\k{\kappa} \def\l{\lambda}
\def\p{\psi}  \def\s{\sigma} \def\ta{\theta} \def\t{\tau}
\def\om{\omega}   \def\DD{\Delta}

\def\bC{{\mathbf C}} \def\bN{{\mathbf N}} 
\def\bR{{\mathbf R}}  \def\bT{{\mathbf T}} \def\bZ{{\mathbf Z}}

 \def\id{\mathrm{id}}
 \def\tr{\mathrm{tr}}
\def\Re{\mathrm{Re\,}}

\def\ip#1#2{{(\,#1\mid #2\,)}}
\def\m#1{{\vert #1 \vert}} \def\n#1{\Vert #1 \Vert}

\def\demo{\noindent \emph{Proof.\ }}

\def\sa{{\cal A}}  \def\AA{(A,\Delta)} 
\def\AAu{(\au,\Delta)}
\def\ar{{A_{{\rm r}}}} \def\hr{{h_{{\rm r}}}} \def\hu{{h_{{\rm u}}}}
\def\er{\e_{{\rm r}}} 
\def\omu{{(\om\otimes \id)(U)}} \def\tu{{(\id\otimes \t)(U)}}
\def\nuu#1{\n{#1}_{{\rm u}}} \def\nr#1{\n{#1}_{{\rm r}}}
\def\nc#1{\n{#1}_{{\rm c}}} \def\np#1{\n{#1}_{{\rm \pi }}}
\def\au{{A_{{\rm u}}}} \def\eu{\e_{{\rm u}}}
\def\rom#1{{\rm #1}} \def\arc{A_{\rom{rc}}} \def\aarc{\A_{\rom{rc}}}
\def\ddrc{\DD_{\rom{rc}}} \def\ntwo#1{\n{#1}_2}
\def\hrc{h_{\rom{rc}}} \def\B{{\cal B}}

\def\A{{\cal A}}

\title{Co-Amenability of Compact Quantum Groups}

\author{E.~B\'edos$^*$\ \   G.J.~Murphy  \ \ L.~Tuset$^*$}

\date{October 20, 2000}

\begin{document}
\maketitle

\begin{abstract}
We study the concept of co-amenability for a compact quantum
group. Several conditions are
derived that are shown to be equivalent to it. Some consequences
of co-amenability
that we obtain are faithfulness of the Haar integral and
automatic norm-boundedness of positive linear functionals on
the quantum group's Hopf $*$-algebra (neither of these
properties necessarily holds without co-amenability).

\vspace{5ex}

\hspace{-\parindent} \emph{Subj. Class.}: Ouantum groups, C*-algebras

\hspace{-\parindent} \emph{MSC 2000}:  Primary 46L05, 46L65.
Secondary 16W30, 22D25, 58B32.

\hspace{-\parindent} \emph{Keywords}: compact quantum group, amenability.
\end{abstract}

\vfill
\thanks{
$^*$  partially supported by the Norwegian Research Council}

\newpage

\section{Introduction}

In this paper we introduce and study a concept of co-amenability for
compact quantum groups defined in the sense of
S.L.~Woronowicz~\cite{Wo1, Wo2}---see also~\cite{MT} for an
exposition that provides much of the background for this
paper.
Co-amenability of so-called regular multiplicative unitaries has
been introduced by S. Baaj and G. Skandalis \cite[Appendix]{BS} \cite{Bla}. 
One can then procede to define co-amenability of a compact
quantum group by requiring that the regular multiplicative unitary
associated to its reduced quantum group is co-amenable. However, the
C*-algebra formulation of compact quantum groups is
more accessible than the theory of multiplicative
unitaries, which is technically quite involved.
We therefore feel that it is 
worthwhile and appropriate to present a direct definition of co-amenability, 
which is perhaps more intrinsic to the $C^*$-algebra theory of 
compact quantum groups. The Baaj-Skandalis approach to co-amenability 
for compact quantum groups has been rephrased by T. Banica \cite{Ba,Ba2}
to accomodate this, but details are defered to Baaj-Skandalis' work. Our
exposition starts from an elementary remark of Woronowicz \cite[p.~623]{Wo1} 
and is aimed to be self-contained. To motivate our 
definition we briefly discuss here the concept of amenability for a 
discrete group and its equivalent formulations in terms of the group
C*-algebras \cite {Pat,Ped}.

If $\Gamma$ is a discrete group, its reduced and full group
C*-algebras $C^*_{{\rm r}}(\Gamma)$ and $C^*(\Gamma)$ can be
endowed with co-multiplications $\DD_{{\rm r}}$ and $\DD$
making them into compact quantum groups. Details are given in
Section~2. We shall call these the {\em reduced} and {\em
universal} compact quantum groups associated with~$\Gamma$.
The Haar integrals of ${(C^*_{{\rm r}}(\Gamma),\DD_{{\rm
r}})}$ and ${(C^*(\Gamma),\DD)}$ are the canonical tracial
states. Since the left kernel of the trace on $C^*(\Gamma)$ is
the kernel of the canonical  \\ $*$-homomorphism $\ta$ from
$C^*(\Gamma)$ onto $C^*_{{\rm r}}(\Gamma)$, faithfulness of
the Haar integral of ${(C^*(\Gamma),\DD)}$ is equivalent to
amenability of $\Gamma$. Of course, we are using here the well
known equivalence of amenability of $\Gamma$ and injectivity
of~$\ta$; this result is often
called the Hulanicki-Reiter theorem in the literature. The
co-unit of ${(C^*(\Gamma),\DD)}$ is norm-bounded, but that of
${(C^*_{{\rm r}}(\Gamma),\DD_{{\rm r}})}$ may not be. In fact,
it is known that $\Gamma$ is amenable if, and only if, the
co-unit of the latter is norm-bounded. This is essentially a
reformulation of the classical result that $\Gamma$ is
amenable if, and only if, the trivial 1-dimensional
representation of $\Gamma$ is weakly contained in the regular
representation.

This discussion serves to motivate our introduction of the
concept of co-amenability for a general compact quantum group
and we shall frequently refer back to these examples for the
purposes of illustration and motivation of the results we
obtain in the sequel. We define a compact quantum
group $\AA$ to be {\em co-amenable} if the co-unit of its reduced
compact quantum group $(\ar,\DD_{{\rm r}})$ is norm-bounded (see
Section~2
for the definition of $(\ar,\DD_{{\rm r}})$).

If a concept is to be a fruitful one in an abstract theory, it
is desirable that it have a number of different formulations.
Indeed we show that co-amenability is equivalent to several
other conditions; one of these equivalences is an
analog of the Hulanicki--Reiter theorem
(see~Theorem~\ref{am}), which etablishes the link with Banica's
definition. One particularly nice condition
ensuring co-amenability of a compact quantum group is the existence
of a non-zero
multiplicative linear functional on its
reduced quantum group
(Corollary~\ref{thm: character corollary}).

A co-amenable compact quantum group has a number of desirable
properties not possessed by arbitrary compact quantum groups.
We show, for example, that a co-amenable compact quantum group
has a faithful Haar integral (it then follows that the Haar
integral is a KMS~state~\cite{KT,KV}). If a compact quantum group
is not co-amenable, then the co-unit on the Hopf $*$-algebra
of its reduced compact quantum group provides an example of a
positive linear functional that is {\em not} norm-bounded.
However, we show that every positive linear functional on the
Hopf $*$-algebras of a co-amenable compact quantum group is
necessarily norm-bounded (Corollary~\ref{pos}).

The use of the word {\em co-amenability}
deserves some explanation. First recall that amenability of
Kac algebras \cite{ES} is defined in terms of the existence
of an invariant state. If we define amenability of a compact
quantum group in these terms, namely by requiring only the
existence of an invariant state, then all compact quantum
groups are trivially amenable, since the Haar integral is an
invariant state. Thus, this is not a satisfactory definition.
On the other hand, the natural concept of amenability for
discrete quantum groups makes good sense---we
study this notion in a forthcoming paper~\cite{BMT}.
There is a relationship between co-amenability of a compact
quantum group as defined in this paper and amenability of the
associated dual discrete quantum group. The chosen terminology
is aimed to reflect this dual relationship. It also fits with the one
introduced by Baaj and Skandalis in \cite{BS} for regular
multiplicative unitaries. Note however the slightly confusing fact
that Banica \cite{Ba, Ba2} uses most of the time the word amenability
instead of
co-amenability for compact quantum groups (which he calls ``Woronowicz
algebras'').

The paper is organized as follows: In Section~2 we construct
the reduced quantum group corresponding to a compact quantum
group and use it to define co-amenability of the original
compact quantum group. We then derive conditions equivalent to
co-amenability and show it
implies faithfulness of the Haar integral.
As an application of the ideas in this section,
we give a new proof of the theorem of G.~Nagy on faithfulness
of the Haar measure of quantum $SU(2)$. In Section~3 we
consider the universal compact quantum group associated to a
compact quantum group and obtain other conditions equivalent
to co-amenability; in addition, we prove the norm-boundedness
result for positive linear functionals alluded to above. Our
final section, Section~4, is a short one in which we explore
the idea of a bounded co-unit in the context of a compact
quantum semigroup and show that if the latter admits a
faithful Haar integral and a bounded co-unit, it is
necessarily a co-amenable compact quantum group.

For the ease of the reader, our account is quite detailed and
we provide  proofs of several important results which are presented
in a rather sketchy manner in the literature. Especially, we give
in an appendix a proof of the uniqueness property of the
associated dense Hopf $*$-algebra of a compact quantum group.
This useful property is stated without proof in \cite{KT}.

We shall use the convention that ${X\otimes Y}$ represents the
algebraic tensor product when $X$ and $Y$ are simply linear
spaces, or $*$-algebras that are not C*-algebras; if $X$ and
$Y$ are Hilbert spaces, ${X\otimes Y}$ represents the Hilbert
space tensor product and if $X$ and $Y$ are C*-algebras,
${X\otimes Y}$ represents the spatial C*-tensor
product~\cite[Chapter~6]{GM}.

\section{The Reduced Quantum Group}
\label{red}

Throughout this section $(A, \Delta )$ denotes a compact
quantum group. Its Haar integral is denoted by~$h$. The
associated  Hopf $*$-algebra is denoted by $\A$,  
the co-inverse by $\kappa$ and the co-unit
by~$\varepsilon$.  Recall that $\e$ and $\k$ are,
in general, only defined on $\A$. One can describe $\A$ by
saying it is the unique  Hopf $*$-algebra for which $\A$
is a dense unital $*$-subalgebra of~$A$ and the co-multiplication
of~$\A$ is
obtained by restriction of the co-multiplication of~$A$.
The reader may find some basic
definitions and a proof of this
uniqueness property in an appendix to this paper.
We refer otherwise to \cite{MT} and
\cite{Wo2} for
the basic theory of compact quantum groups.

Let $(C(G),\DD)$ be  a commutative compact
quantum group associated to a compact group $G$, the co-multiplication
$\DD$ being  dual to the group multiplication operation ${G\times
G\to G}$. In this case the Haar integral $h$ is the integral with respect to
the
Haar measure on $G$. This has full support and therefore $h$
is faithful. Faithfulness of the Haar integral no longer holds
for an arbitrary compact quantum group. To illustrate this
we return to the group C*-algebras
of a discrete group and discuss them in a little more detail.

Let $\Gamma$ be a discrete group and let ${L\colon x\mapsto
L_x}$ be the left regular representation of $\Gamma$
on~$\ell^2(\Gamma)$. Thus, if $(\d_x)_{x\in \Gamma}$ is the
canonical orthonormal basis of $\ell^2(\Gamma)$,
$L_x(\d_y)=\d_{xy}$. Let $\ar=C^*_{\rm r}(\Gamma)$ be the
reduced group C*-algebra of $\Gamma$; that is, $\ar$ is the
C*-subalgebra of $B(\ell^2(\Gamma))$ generated by the
operators $L_x$ (${x\in \Gamma}$). The linear map $\DD_{{\rm r}}$
defined on~$\ar$ by $\DD_{{\rm r}}(L_x)={L_x\otimes L_x}$, for all
$x\in
\Gamma$, is a co-multiplication of $\ar$. (To see that $\DD_{{\rm r}}$
is well defined, observe that there is a unitary operator $W$
on ${\ell^2(\Gamma)\otimes \ell^2(\Gamma)}$ for which
${L_x\otimes L_x}={W^*(1\otimes L_x)W}$, for all $x\in
\Gamma$; $W$ is defined by setting ${W(\d_x\otimes
\d_y)}={\d_{x^{-1}y}\otimes \d_y}$, for all $x,y\in \Gamma$.)
It is easy to see that ${(\ar\otimes 1)\DD_{{\rm r}} \ar}$ and
${(1\otimes \ar)\DD_{{\rm r}} \ar}$ each have closed linear span equal
to ${\ar\otimes \ar}$. Hence, $(\ar,\DD_{{\rm r}})$ is a compact
quantum
group.

It is well known that $C^*_{\rm r}(\Gamma)$ admits a faithful
tracial state $\tr$ given by ${\tr(L_x)=0}$, if $x$ is an
element of $\Gamma$ that is not equal to the unit of $\Gamma$.
In fact, $\tr$ is the Haar integral of
$(\ar,\DD_{{\rm r}})$~\cite[Example~10.4]{MT}. The  dense Hopf $*$-algebra
$\A_{{\rm
r}}$ of $(\ar,\DD_{{\rm r}})$ is the linear span of all the unitaries
$L_x$
($x\in \Gamma$). It may be identified with the group algebra

$\bf{C}$($\Gamma$)
of $\Gamma$ equipped with its canonical Hopf
$*$-algebra structure.

The full group C*-algebra $\au=C^*(\Gamma)$ is, by definition,
the enveloping C*-algebra of the Banach $*$-algebra $
\ell^{1}( \Gamma)$. By construction, $\bf{C}$($\Gamma$) 
is dense
in~$\au$. Therefore, $\Gamma$ admits a universal unitary
representation, ${V\colon \Gamma\to \au}$, ${x\mapsto V_x}$
such that the linear span of the elements~$V_x$ is dense
in~$\au$. A co-multiplication on $\au$ making it into a
compact quantum group is determined by first setting
$\DD(V_x)=V_x\otimes V_x$, for all $x\in \Gamma$, and then
extending $\DD$ to $\au$ by its universal property. The Hopf
$*$-algebra $\A_{{\rm u}}$ of $\AAu$ is the linear span of the
elements~$V_x$, and it too may be identified with
$\bf{C}$($\Gamma$).

By the universal property of $C^*(\Gamma)$ there exists a
canonical surjective $*$-homomorphism $\theta : \au
\rightarrow \ar $ mapping each $V_{x}$ onto $L_{x}$, hence
mapping  $\A_{{\rm u}}$ onto $\A_{{\rm r}}$.  The Haar integral
on $\au$ is the canonical tracial state of $\au$ given by $h =
\tr \circ\theta$. Its left kernel $N_h$ is clearly the kernel
of $\ta$,  so $\ar=\au/N_h$. Again using the universal
property of $C^*(\Gamma)$, we see there is a
$*$-homomorphism~$\e$ from $\au$ to $\bC$ such that
$\e(V_x)=1$, for all $x\in \Gamma$. A simple computation shows
that $\e$ is the co-unit for $\AAu$. (More precisely, the
restriction of $\e$ to the Hopf $*$-algebra of $\AAu$ is the
co-unit.) The important point here is that $\e$ is
norm-bounded.

The group $\Gamma$ is amenable if, and only if, $\theta$ is
injective, and the co-unit of $C^*_{\rm r}(\Gamma)$ is
therefore norm-bounded in this case. If $\Gamma$ is not
amenable, this co-unit is not norm-bounded, as pointed out in
the Introduction. In the case that $\Gamma={\bf F}_2$, the
free group on two generators, one can see the co-unit of
$C^*_{{\rom r}}(\Gamma)$ is not norm-bounded by means of the
well known fact that $C^*_{\ r}(\Gamma)$ is simple (and not
one-dimensional!) and therefore admits no $*$-homomorphism
onto $\bC$.

\vspace{2ex}
Suppose now that $\AA$ is an arbitrary compact quantum group
with associated Hopf $*$-algebra $\A$. It is known \cite{Wo2} that the
Haar integral of $(A, \DD )$ is faithful on $\A$, but as we have seen,
in general, not on the $C^*$-algebra $A$. We will now furnish a 
$C^*$-algebra envelope of the Hopf $*$-algebra $\A$ for which the
Haar integral is faithful.
Recall that the left kernel $N_h$ of~$h$ is a two-sided ideal
of $A$~\cite{Wo2}. Set $\ar=A/N_h$ and let $\theta$ be the
quotient map from $A$ onto $\ar$. We shall make $\ar$ into a
compact quantum group. This reduction procedure is sketched 
in \cite{Wo1}, but no details are given there, or anywhere 
else in the literature that we are aware of. Since this is an 
important construction for this paper
we give the required details in the following result.

\begin{thm}
If $\AA$ is a compact quantum group, then the C*-algebra $\ar$
can be made into a compact quantum group whose
co-multiplication $\DD_{{\rm r}}$ is determined by $\DD_{{\rm
r}} (\ta(a))=(\theta\otimes \theta)\DD(a)$, for all $a\in A$.
The Haar
integral of $(\ar,\DD_{{\rm r}})$  is the unique state
$\hr$ of $\ar$ such that $h=\hr\circ\ta$. The state $\hr$ is
faithful. Also, the quotient map $\ta$ is faithful on~$\A$ and 
the Hopf $*$-algebra of $(\ar,\DD_{{\rm r}})$ is
$\ta(\A)$, with co-unit $\er$ and co-inverse $\kappa_{{\rm r}}$ 
determined by $\e = \er \circ \ta$ and $ \ta \circ \kappa = \kappa_{{\rm r}} 
\circ \ta$, respectively.
\end{thm}

\demo To show that we can define a $*$-homomorphism
${\DD_{{\rm r}}\colon \ar\to \ar\otimes \ar}$ such that
$\DD_{{\rm r}}(\ta(a))=(\theta\otimes \theta)\DD(a)$, for all
$a\in A$, we need only show that $\ker(\ta)\subseteq
{\ker(\ta\otimes\ta)\DD)}$. Clearly, it suffices to show that
$\ker(\ta)\subseteq {\ker((\id\otimes\ta)\DD)}$. To see this,
we first observe that, by the Cauchy-Schwartz inequality, $h$ vanishes
on~$\ker(\ta)$. Therefore it induces a unique state $\hr$ on~$\ar$
such that $h=\hr \circ \ta$. Since $\ker(\ta)=N_h$, it is
clear that $\hr$ is faithful. Using the fact that product states
separate elements of $\ar\otimes \ar$, it easily follows that
${\id\otimes \hr\colon \ar\otimes \ar\to \ar}$ is faithful.
Suppose now $\ta(a)=0$.
Then $h(a^*a)=0$ and therefore, ${(\id\otimes \hr)(\id\otimes
\ta)\DD(a^*a)} \\ ={(\id\otimes h)\DD(a^*a)}={h(a^*a)1=0}$.
Consequently, ${(\id\otimes \ta)\DD(a^*a)=0}$, and therefore
${(\id\otimes \ta)\DD(a)~=~0}$ as required. Thus, we can well
define a $*$-homomorphism $\DD_{{\rm r}}$ as claimed above.

One can easily check now that $\DD_{{\rm r}}$ is a
co-multiplication on~$\ar$. Since the linear spans of
${(1\otimes A)\DD(A)}$ and ${(A\otimes 1)\DD(A)}$ are dense in
${A\otimes A}$, it follows immediately that the linear spans
of ${(1\otimes \ar)\DD_{{\rm r}}(\ar)}$ and ${(\ar\otimes
1)\DD_{{\rm r}}(\ar)}$ are dense in ${\ar\otimes \ar}$. Hence,
$(\ar,\DD_{{\rm r}})$ is a compact quantum group.

If $a\in A$, then ${(\id \otimes \hr)\DD_{{\rm r}}(\ta(a))}=
{(\id\otimes \hr)(\ta\otimes \ta)\DD(a)}={\ta(\id\otimes
h)\DD(a)}={\ta(h(a)1)}=\hr(\ta(a))\ta(1)$. Similarly, ${(\hr
\otimes \id)\DD_{{\rm r}}(\ta(a))}=\hr(\ta(a))\ta(1)$. Hence,
$\hr$ is the Haar integral of $(\ar,\DD_{{\rm r}})$.

The injectivity of $\ta$ on~$\A$ follows readily: If $a\in \A$
and $\ta(a)=0$, then ${h(a^*a)=0}$. Since $h$ is faithful
on~$\A$, we deduce that $a=0$.

We can therefore define linear maps, $\er : \ta(\A) \rightarrow \bC$ and 
$\kappa_{{\rm r}} : \ta(\A) \rightarrow \ta(\A)$, by setting  
$ \er (\ta(a)) = \e(a)$ and $ \kappa_{{\rm r}} (\ta(a)) = \ta(\kappa(a))$, 
for all $a \in \A $.  
It is then clear that $\ta(\A)$ is a dense Hopf
$*$-subalgebra of $(\ar,\DD_{{\rm r}})$ with co-unit $\er$ and 
co-inverse $\kappa_{{\rm r}}$.
Hence, by uniqueness,
$\ta(\A)$ is the Hopf $*$-algebra associated
to $(\ar,\DD_{{\rm r}})$.~\qed

We call the compact quantum group $(\ar,\DD_{{\rm r}})$ described in the
theorem  the {\em reduced quantum group}
of~$(A,\Delta)$ and we call $\ta$ the {\em canonical map} from
$A$ onto~$\ar$. It is clear that $\theta$ is a $*$-isomorphism if, and only
if, $h$ is faithful.

\vspace{2ex}
If $\AA$ is the universal compact quantum group associated to
a discrete group~$\Gamma$, then the reduced compact quantum
group of~$\AA$ is equal to the reduced compact quantum group
of~$\Gamma$; that is, $(\ar,\DD_{{\rm r}})=(C^*_{{\rom
r}}(\Gamma),\DD_{{\rom r}})$. That $\ar=C^*_{{\rm r}}(\Gamma)$
follows from the fact that the left kernel of the Haar
integral of $\AA$ is equal to the kernel of the canonical
$*$-homomorphism $\ta$ from $C^*(\Gamma)$ onto $C^*_{{\rm
r}}(\Gamma)$, as we have observed before. The only other item
that needs to be checked is that $\DD_{{\rom
r}}\ta={(\ta\otimes \ta)\DD}$, and this easily follows from
the definitions of the co-multiplications on $C^*(\Gamma)$ and
$C^*_{{\rm r}}(\Gamma)$.

\vspace{2ex}

If $\AA$ is an arbitrary compact quantum group, we say it is
{\em co-amenable} if the co-unit $\er$ of $(\ar,\DD_{{\rm
r}})$ is norm-bounded. We can then extend the co-unit to a  
$*$-homomorphism $\er$ on~$\ar$. A consequence is that $A$
is never simple, if $\AA$ is co-amenable, since the kernel of
$\er\ta$ is a closed two-sided ideal of $A$ of
co-dimension one.

From our discussion above, it is evident that the reduced (resp.
universal) compact quantum group associated to a discrete
group $\Gamma$ is co-amenable if, and only if, $\Gamma$ is
amenable. Note also that a finite quantum group---that is, a compact
quantum group $\AA$ for which $A$ is finite dimensional---is
necessarily co-amenable, since in this case $A=\A$.

It is perhaps of some interest to interpret the idea of
co-amenability in the context of a commutative compact quantum group
$(C(G),\DD)$ associated to a classical compact
group~$G$. Since the Haar integral is faithful, as we observed
before, $(C(G),\DD)$ is co-amenable if its co-unit is
norm-bounded. That this is the case is trivial, since the
co-unit is given by the (restriction of) the
evaluation map, ${f\mapsto f(e)}$, where $e$ is the unit
of~$G$. Thus, a classical compact group is ``co-amenable''.

The following theorem allows us to verify co-amenability
without reference to the reduced compact quantum group.
However, its real importance is its assertion that
faithfulness of the Haar integral is a consequence of
co-amenability. In practice, it provides a useful method of
showing such faithfulness (see Corollary~\ref{thm; Nagy's
theorem} below).

The first paragraph of the proof of the theorem is taken from
the proof of Theorem~8.1 of~\cite{MT} (the exactness
assumption on~$A$ used in~\cite{MT} is not needed here).

\begin{thm} \label{thm: Haar faithfulness} A compact
quantum group $\AA$ is co-amenable if, and only if, its Haar
integral is faithful and its co-unit is norm-bounded.
\end{thm}

\demo Clearly, we need only show that if $\AA$ is co-amenable,
then $h$ is faithful. Let $I=N_h$. If $a\in I$ and $\s$ is a
positive linear functional on~$A$, then ${(\s\otimes
h)\DD(a^*a)}=\s(1)h(a^*a)=0$, since ${(\id\otimes
h)\DD(a^*a)=h(a^*a)1}$. Hence, since ${\s\otimes h}$ is
positive, ${(\s\otimes h)(c\DD(a))=0}$, for all ${c\in
A\otimes A}$. Because $\s$ is an arbitrary positive linear
functional on~$A$, this implies  ${(\id\otimes
h)(c\DD(a))=0}$. If $\t\in A^*$ and $c={1\otimes b}$, where
$b\in A$, then we have ${h(b(\t\otimes \id)(\DD
a))}={\t((\id\otimes h)(c\DD(a)))=0}$. Hence,
$(\t\otimes\id)\DD(a)\in I$.

The co-units $\er$ and $\e$ are norm-bounded, by
 co-amenability, so admit extensions $\er$ and $\e$ to $\ar$ and $A$,
respectively, which satisfy
$\e=\er\ta$. It follows that $\t(a)={\t((\id\otimes
\e)\DD(a))}={\er\theta((\t\otimes
\id)\DD(a))}={\er(0)=0}$. Since $\t$ was an arbitrary
element of $A^*$, we must have $a=0$. Hence, $N_h=I=0$; that
is, $h$ is faithful.~\qed

It follows from Theorem~\ref{thm: Haar faithfulness} that co-amenability is
preserved under formation of the tensor product of two compact
quantum groups. This is the quantum counterpart of the
statement that a product of two discrete amenable groups is
amenable. Recall that the {\em tensor product} of two
compact quantum groups ${(A_i , \Delta_i )}$ is the compact
quantum group $(A,\DD)={(A_1\otimes A_2,
\Delta_1\times\Delta_2 )}$ with comultiplication defined by
\[\Delta_1\times\Delta_2 =(\id\otimes F\otimes
\id)(\Delta_1\otimes\Delta_2 ):A\rightarrow A\otimes A,\]
where $F:A_1\otimes A_2\rightarrow A_2\otimes A_1$ denotes the
flip map given by ${F(a_1\otimes a_2)} = {a_2\otimes a_1}$,
for $a_1\in A_1$ and $a_2\in A_2$. The Hopf $*$-algebra of
$\AA$ is $\A_1\otimes \A_2$, where $\A_i$ is the Hopf
$*$-algebra of~$(A_i,\DD)$; the Haar integral and the co-unit
of $\AA$ are $h_1\otimes h_2$ and
$\varepsilon_1\otimes\varepsilon_2$, respectively, where $h_i$
is the Haar integral and $\varepsilon_i$ is the co-unit of
$(A_i ,\Delta_i )$.

If $(A_i , \Delta_i )$ are both co-amenable, then, by Theorem
\ref{thm: Haar faithfulness}, their Haar integrals $h_i$ are
faithful, and therefore $h_1\otimes h_2$ is also faithful.
Hence $\AA$ is equal to its reduced compact quantum group, so
we only need to check that the co-unit
$\varepsilon_1\otimes\varepsilon_2$ is norm-bounded and this
is obvious, since $\varepsilon_i$ are both norm-bounded. Thus,
$\AA$ is co-amenable.

In the reverse direction, if $\AA$ is co-amenable, then both
$(A_1,\DD_1)$ and $(A_2,\DD_2)$ are co-amenable. For,
faithfulness of ${h_1\otimes h_2}$ trivially implies
faithfulness of each of $h_1$ and $h_2$; equally easily,
norm-boundedness of ${\e_1\otimes \e_2}$ implies
norm-boundedness of $\e_1$ and $\e_2$. Hence, co-amenability
of $(A_1,\DD_1)$ and $(A_2,\DD_2)$ follows from
Theorem~\ref{thm: Haar faithfulness}.

This observation allows
us to give an example of a compact quantum group $\AA$ that is
not co-amenable and that is neither co-commutative nor
commuta-\\ tive: We set $A_1=C^*({\bf F}_2)$ and $A_2=C({\bf
S}_3)$, where ${\bf F}_2$ is the free group on two generators 
and ${\bf S}_3$ is the finite (compact) group of
permutations on three symbols. Then  we let  $\AA$ be the
tensor product of these two compact quantum groups.

\vspace{2ex}

We turn now to finding other conditions equivalent to
co-amenability or, more generally, conditions equivalent to
norm-boundedness of the co-unit $\e$.

Recall  a finite-dimensional unitary co-representation
${U\in M_N(\bC)\otimes \A}$ of $\AA$ is said to be {\em
fundamental} if its matrix elements $U_{ij}$ (relative to some
system of matrix units for $M_N(\bC)$) generate the Hopf
$*$-algebra $\A$ associated to $\AA$, as a $*$-algebra. The
{\em compact matrix pseudogroups}, as defined by Woronowicz
in \cite{Wo1},
are precisely the compact quantum groups that admit a
fundamental unitary co-representation.

The equivalence of Conditions~(1) and~(2) in the corollary of
the following theorem can be regarded as a generalization of
H.~Kesten's classical characterization of the amenability of a
finitely-generated discrete group in terms of the spectrum of
the sum of the generators in the regular representation (see
\cite{Kes}, and also \cite{HRV}). This equivalence, which is
due to G.~Skandalis, is proved in \cite{Ba}.
Its connection to Kesten's result is explained in~\cite{Ba2}.
The proof of our more general result is somewhat different.

\begin{thm} \label{thm: boundedness theorem}
Suppose that $\AA$ is a compact matrix pseudogroup and that
${U\in M_N(\bC)\otimes \A}$ is a fundamental unitary
co-representation of~$\AA$.\\ We set $\chi_U=\sum_{i=1}^N U_{ii}$.\\ 
Of course, since $\n{U_{ij}}\le 1$, for all indices $i$ and
$j$, $\n{\Re\chi_U}\le N$.\\ 
The following are equivalent conditions:

\vspace{2ex}
(1) The co-unit $\e$ of $\AA$ is norm-bounded;

(2) $N$ belongs to the spectrum of $\Re\chi_U$ in $A$;

(3) There exists a state $\t$ on $A$ such that $\t(\Re
\chi_U)=N$;

(4) There exists a state $\t$ on $A$ such that $\t(U_{ii})=1$,
for ${i=1,\dots,N}$.

(5) For all scalars ${\l_0,\l_1,\dots,\l_N}$, \begin{equation}
\label{eqn: norm inequality} \m{\sum_{i=0}^N\l_i}\le \n{\l_0
1+\sum_{i=1}^N \l_iU_{ii}}.\end{equation}

\end{thm}

\demo Recall first from~\cite[Proposition~1.8]{Wo1} that $\e$ is
uniquely determined on $\A$ by $\e(U_{ij})=\delta_{ij}$, for
all indices $i$ and $j$. Especially, $\e(U_{ii})=1$ for all $i$, 
so  we have
$\sum_{i=0}^N\l_i=\e(\l_{0}1 + \sum_{i=1}^N \l_{i}U_{ii})$.
The implication $(1)\implies (5)$ follows by noting that if
$\e$ is norm-bounded, its norm must be equal to one, and
Inequality~(\ref{eqn: norm inequality}) is an immediate
consequence. To see Condition~(5) implies~(4), we note that
Inequality~(\ref{eqn: norm inequality}) implies that the
linear functional  $\t_0$, defined on the linear span of 1 and
the elements $U_{ii}$ by mapping all of these elements to 1 in
$\bC$, is well defined and has norm equal to~1. By the
Hahn--Banach theorem, $\t_0$ extends to a norm-one linear
functional~$\t$ on $A$. Since $\t(1)=\n\t=1$, $\t$ is a state
of $A$.

Since a state is necessarily self-adjoint, the implication
$(4)\implies (3)$ is clear.

 Set $X_{ij}=U_{ij}-\d_{ij}$ and
${X=\sum_{i,j=1}^N X_{ij}^*X_{ij}+X_{ij}X_{ij}^*}$. Using the
fact that $\sum_{i=1}^N U^*_{ij}U_{ij}=\sum_{i=1}^N
U_{ij}U_{ij}^*=1$, we have ${X= 4(N-\Re\chi_U)}$. Hence, the
element ${N-\Re\chi_U}$ is positive. Therefore,
${N-\Re\chi_U}$ is invertible if, and only if, there exists a
positive number such that $N-\Re \chi_U\ge \d$. Hence, $N$
belongs to the spectrum of $\Re\chi_U$ if, and only if,
$\t(\Re\chi_U)=N$, for some state $\t$ of $A$. That is,
Conditions~(2) and~(3) are equivalent.

Thus, it remains only to show that $(3)\implies (1)$. Suppose
Condition~(3) holds, so that there exists a state $\t$ on $A$
such that $\t(N-\Re\chi_U)=0$ and therefore, $\t(X)=0$. Hence,
${\t(X^*_{ij}X_{ij})=\t(X_{ij}X^*_{ij})=0}$. Let $\f$ be the
GNS representation associated to $\t$, acting on the Hilbert
space $H$, and let $x$ be the canonical cyclic vector
associated to this representation, so that
$\t(a)=\ip{\f(a)x}x$ and $\f(A)x$ is dense in $H$. Clearly,
${\f(X_{ij})x=\f(X^*_{ij})x=0}$ and therefore
$\f(U_{ij})x=\f(U^*_{ij})x=\d_{ij}x$. Hence, if $a$ is product
of matrix elements $U_{ij}$ and $U^*_{kl}$, then $\f(a)x\in
\bC x$. Since $U$ is a fundamental co-representation of~$\AA$,
the closed linear span of such products is equal to~$A$ and
therefore $\f(A)x\subseteq \bC x$. Hence, $H=\bC x$ and
therefore $\dim(H)=1$. It follows that $\f$ is scalar-valued
and therefore $\f(a)=\t(a)1$, for all $a\in A$. Hence, $\t$ is
a norm-bounded $*$-homomorphism. Moreover, since
$\m{\t(X_{ij})}^2\le \t(X^*_{ij}X_{ij})=0$, we have
$\t(U_{ij})=\d_{ij}=\e(U_{ij})$, for ${i,j=1,\dots,N}$. Hence,
since the elements $U_{ij}$ generate $\sa$ as a $*$-algebra,
$\t=\e$ on $\sa$ and therefore $\e$ is norm-bounded.~\qed

\begin{cor} \label{thm: co-amenability corollary} With the
same assumptions as in the preceding theorem, the following
are equivalent conditions:

(1) $\AA$ is co-amenable;

(2) $N$ belongs to the spectrum of $\ta(\Re\chi_U)$ in $\ar$;

(3) There exists a state $\t$ on $\ar$ such that $\t\ta(\Re
\chi_U)=N$;

(4) There exists a state $\t$ on $\ar$ such that
$\t\ta(U_{ii})=1$, for ${i=1,\dots,N}$.

(5) For all scalars ${\l_0,\l_1,\dots,\l_N}$, \[
 \m{\sum_{i=0}^N\l_i}\le \n{\l_0
1+\sum_{i=1}^N \l_i\ta(U_{ii})}.\]

\end{cor}

\demo The result follows from the theorem by observing that
$(\id\otimes\ta)(U)$ is a fundamental co-representation of
$(\ar,\DD_{{\rm r}})$.~\qed

If $U$ is a a unitary co-representation of~$\AA$ on a Hilbert
space~$H$, so that $U\in M(K(H)\otimes A)$, the multiplier
algebra of ${K(H)\otimes A}$, recall that its matrix elements
are the elements of~$A$ of the form~$\omu$, where $\om$ is a
strictly continuous linear map on $K(H)$. Not every compact
quantum group admits a fundamental unitary co-representation
but all admit a unitary co-representation for which the matrix
elements generate its C*-algebra (for example, the matrix
elements of the  regular co-representation have dense
linear span in the C*-algebra).

If $U$ is any unitary co-representation of~$\AA$ on a Hilbert
space $H$ and the co-unit $\e$ is norm-bounded, then
${(\id\otimes\e)(U)=1}$ in $B(H)$. For, the equality
${(\id\otimes \e)\DD=\id}$ implies
${U=(\id\otimes (\id\otimes \e)\DD)(U)}={(\id\otimes\id\otimes
\e)(\id\otimes \DD)(U)}={(\id\otimes \id\otimes
\e)(U_{(12)}U_{(13)})}= {U((\id\otimes\e)(U)\otimes 1)}$.
Since $U$ is invertible, we deduce that ${1\otimes
1=(\id\otimes\e)(U)\otimes 1}$ and therefore
${(\id\otimes\e)(U)= 1}$, as required.

The following represents a partial generalization
of~Theorem~\ref{thm: boundedness theorem}.

\begin{thm} Let $U$ be a unitary co-representation of $\AA$ whose
matrix elements generate $A$ as a C*-algebra. Then the
following are equivalent conditions:

(1) The co-unit $\e$ is norm-bounded;

(2) There exists a state $\t$ of $A$ for which $\tu=1$.
\end{thm}

\demo Taking $\t=\e$, the implication $(1)\implies (2)$ is
immediate from the remarks preceding this theorem. To see the
converse, suppose given a state $\t$ of $A$ for which $\tu=1$.
Let $\f$ be the GNS representation associated to $\t$. We
suppose $H$ is the Hilbert space on which $\f$ acts and that
$x$ is the canonical cyclic vector for~$\f$. As in the proof
of Theorem~\ref{thm: boundedness theorem}, we shall show that
$\f(a)=\t(a)1$, for all $a\in A$. First, let $a=\omu$ be a
matrix element of $U$, where $\om$ is a strictly continuous
linear map on $K(H)$. We shall show that $\f(a)x,\f(a)^*x\in
\bC x$. Since $\om$ is linear combination of strictly
continuous states on $K(H)$, to show this result we may
suppose that $\om$ is a state. Then $\n a\le 1$ and
$\t(a)=\om(\tu)=\om(1)=1$; hence, ${0\le
\t((a-1)^*(a-1))}={\t(a^*a)-\t(a)-\t(a)^-+\t(1)}\le
{\t(1)-1-1+\t(1)=0}$. Consequently, ${\t((a-1)^*(a-1))=0}$,
from which it follows that  $\f(a)x=x$. Similar reasoning
shows that ${\t((a-1)(a-1)^*)=0}$ and therefore $\f(a)^*x=x$.
Since the elements $a=\omu$ generate $A$, as a C*-algebra, we
can now argue again as in the proof of Theorem~\ref{thm:
boundedness theorem} to deduce that $\f(A)x=\bC x$. Hence,
$\f(a)=\t(a)1$, for all $a\in A$, as claimed. This implies
that $\t$ is a $*$-homomorphism on~$A$.

Now we shall show that ${(\id\otimes \t)\DD(a)=a}$, for all
$a\in A$. To see this, we may clearly suppose that $a$ is a
matrix element, $a=\omu$, say. Then
\begin{eqnarray*}{(\id\otimes \t)\DD(a)}&=&{(\id\otimes
\t)(\om\otimes
\id\otimes\id)(\id\otimes\DD)(U)}\\&=&{(\id\otimes
\t)(\om\otimes \id\otimes\id)(U_{(12)}U_{(13)})}\\
&=&(\om\otimes \id)(\id\otimes\id\otimes
\t)(U_{(12)}U_{(13)})\\ &=&{(\om\otimes \id)(U((\id\otimes
\t)(U)\otimes 1))}\\&=&{(\om\otimes \id)(U(1\otimes
1))=a}.\end{eqnarray*}

We complete the proof now by showing that $\t(a)=\e(a)$, for
all $a\in \sa$: We have $\t(a)=\t((\e\otimes
\id)(\DD(a))=\e((\id\otimes \t)(\DD(a))=\e(a)$. Hence, $\t$ is a
norm-bounded linear map extending $\e$  and therefore $\e$ is
norm-bounded.~\qed

Let us note explicitly that our proof of the preceding theorem
shows that if $\t$ is as in Condition~(2), then $\t$ is
the---necessarily unique---extension of $\e$ to~$A$.

\begin{cor} Let $U$ be a unitary co-representation of $\AA$ whose
matrix elements generate $A$ as a C*-algebra. Then the
following are equivalent conditions:

(1) $\AA$ is co-amenable;

(2) There exists a state $\t$ of $\ar$ for which ${(\id\otimes
\t\ta)(U) =1}$.
\end{cor}

\demo The element $V=(\id\otimes \ta)(U)$ in the multiplier
algebra  ${M(K(H)\otimes \ar)}$ is a unitary co-representation
of $(\ar,\DD_{{\rm r}})$ whose matrix elements
${(\om\!\otimes\id)(V)}\!=\!{\ta(\omu)}$ generate $\ar$ as a
C*-algebra. The result therefore follows directly from the
theorem.~\qed

We stated before the preceding theorem that it is a partial
generalization of~Theorem~\ref{thm: boundedness theorem}. To
see why, let ${U\in M_N(\bC)\otimes \A}$ be a 
a finite-dimensional unitary co-representation
 of $\AA$ with matrix elements $U_{ij}$ (relative to some
system of matrix units for $M_N(\bC)$).
The equation
$\tu=1$ is clearly equivalent to the condition that
$\t(U_{ij})=\delta_{ij}$, for all $i$ and $j$. Reasoning as in
the proof of Theorem~\ref{thm: boundedness theorem}, this is
easily seen to be equivalent to the condition that
$\t(U_{ii})=1$, for all~$i$. Hence, the preceding theorem
implies the equivalence of Conditions~(1) and~(4)
of~Theorem~\ref{thm: boundedness theorem}.

\vspace{1ex}
We shall need the following result for the proof
of~Theorem~\ref{thm: character theorem}.

\begin{lem} \label{thm: faithfulness lemma}
Let $\AA$ be a compact quantum group for which the Haar
integral~$h$ is faithful. Let $\pi$ be a non-zero
$*$-homomorphism from $A$ to a C*-algebra~$B$. Then the
$*$-homomorphism, ${\hat \pi\colon A\to A\otimes B}$,
${a\mapsto (\id\otimes \pi)\DD(a)}$,  is isometric.
\end{lem}

\demo Let $a\in A$ and suppose that $\hat\pi(a)=0$. Then
$\hat\pi(a^*a)=0$ and therefore ${0=(h\otimes
\id)\hat\pi(a^*a)}={\pi((h\otimes \id)\DD(a^*a))}=
{\pi(h(a^*a)1)}={h(a^*a)\pi(1)}$. Consequently, since
$\pi(1)\ne 0$, we have  $h(a^*a)=0$; faithfulness of~$h$ now
gives $a=0$. Hence, $\hat\pi$ is injective and therefore
isometric.~\qed

The corollary to the following theorem gives another
characterization of co-amenability, this time in terms of a
scalar-valued $*$-homomorphism on the C*-algebra of the
reduced quantum group:

\begin{thm} \label{thm: character theorem}
Let $\AA$ be a compact quantum group for which the Haar
integral~$h$ is faithful.  Then the following are equivalent
conditions:

(1) The co-unit $\e$ is norm-bounded;

(2) There exists a non-zero $*$-homomorphism ${\t\colon A \to
\bC}$.

\end{thm}

\demo The implication $(1)\implies (2)$ is obvious. Suppose
therefore that we have a non-zero $*$-homomorphism ${\t\colon
A\to \bC}$. If $U$ is an $N$-dimensional unitary
co-representation of~$\AA$, then ${(\id\otimes
\t)\DD(U_{ij})}={(\id\otimes \t)(\sum_{k=1}^N U_{ik}\otimes
U_{kj})}={\sum_{k=1}^N U_{ik}\t(U_{kj})}$. Also, since the
matrix $(\t(U_{ij}))$ is a unitary, because $\t$ is a
$*$-homomorphism,  ${\sum_{j=1}^N(\id\otimes
\t)\DD(U_{ij})\t(U_{lj})^-}= {\sum_{k,j=1}^N
U_{ik}\t(U_{kj})\t(U_{lj})^-}={\sum_{k=1}^NU_{ik}\d_{kl}}=U_{il}$.
Hence, recalling that $\sa$ is the linear span of the matrix
elements of finite-dimensional unitary co-representations of
$\AA$, it is clear that the $*$-homomorphism ${\hat\t\colon
\sa\to \sa}$, defined by setting $\hat\t(a)={(\id\otimes
\t)\DD(a)}$, is surjective. Since $h$ is assumed to be
faithful, it follows from Lemma~\ref{thm: faithfulness lemma}
that $\hat\t$ is an isometry. Therefore, if $a\in \sa$,
$\m{\e(\hat\t(a))}= {\m{\t((\e\otimes
\id)\DD(a))}}=\m{\t(a)}\le \n a=\n{\hat\t(a)}$. Therefore,
$\e$ is norm-bounded. Hence, $(2)\implies (1)$.~\qed

\begin{cor} \label{thm: character corollary}
If $\AA$ is an arbitrary compact quantum group, the following
are equivalent conditions:

(1) $\AA$ is co-amenable;

(2) There exists a non-zero $*$-homomorphism ${\t\colon \ar
\to \bC}$.

(3) The Haar integral on $\AA$ is faithful and there exists a
non-zero

\hspace{4 ex} $*$-homomorphism ${\t\colon A \to \bC}$.

\end{cor}

\demo The equivalence between (1) and (2) follows immediately
from the theorem, since the Haar integral of $(\ar,\DD_{{\rm r}})$
is faithful.
The equivalence between (1) and (3) follows by combining the
theorem and Theorem \ref{thm: Haar faithfulness}.~\qed

As an immediate consequence of the equivalence between (1) and
(2) above, we obtain the following corollary which is a
special case of a result in \cite{Be}.

\begin{cor} Let $\Gamma$ be a discrete group.
The following
are equivalent conditions:

(1) $\Gamma$ is amenable;

(2) There exists a non-zero $*$-homomorphism ${\t\colon
C^*_{{\rm r}}(\Gamma)\to\bC}$.

\end{cor}

If $\Gamma$ is a discrete group, then its reduced group
C*-algebra is given by a concrete faithful representation on
the  Hilbert space~$\ell^2(\Gamma)$. Given a compact quantum
group $\AA$, there is a natural faithful
representation of $(\ar,\DD_{{\rm r}})$ whose existence may be
deduced from \cite{BS}. For completeness, we now present this
representation in details.
Let ${\pi\colon A\to B(H)}$ be the
GNS~representation of $A$ associated to the Haar integral~$h$
of $\AA$ and let $z$ be its canonical cyclic vector, so that
$\pi(A)z$ is dense in $H$ and  $h(a)=\ip{\pi(a)z}z$, for all
$a\in A$. We denote by $\ntwo{\cdot}$ the norm of $H$. We set
$\arc=\pi(A)$ and $\aarc=\pi(\A)$, so that $\arc$ is a unital
C*-subalgebra of $B(H)$ and $\aarc$ is a dense unital
$*$-subalgebra of $\arc$. The map~$\pi$ is injective on~$\A$.
For, if $a\in \A$ and $\pi(a)=0$, then
$\ntwo{\pi(a)z}^2=h(a^*a)=0$ and therefore, by faithfulness
of~$h$ on $\A$, $a=0$. Hence, we can define linear maps,
${\ddrc:\aarc\rightarrow\aarc\otimes\aarc}$,
$\e_{{\rm rc}} : \aarc \rightarrow \bC$ and 
$\kappa_{{\rm rc}} : \aarc \rightarrow \aarc$ by setting 
$\ddrc(\pi (a))={(\pi\otimes\pi )\Delta (a)}$,
$\e_{{\rm rc}} (\pi(a)) = \e(a)$ and 
$\kappa_{{\rm rc}} (\pi(a)) = \pi (\kappa(a))$, for all $a \in \A $.  
Clearly, $\ddrc$ is a unital $*$-homomorphism.

\begin{thm} \label{thm: representation theorem}
Let $\AA$ be a compact quantum group and retain the notation
of the preceding paragraph. The map
$\ddrc:\aarc\rightarrow\aarc\otimes\aarc$ has a unique
extension to a $*$-homomorphism ${\ddrc:\arc\rightarrow
\arc\otimes \arc}$. The pair $(\arc ,\ddrc )$ is a compact
quantum group with faithful Haar state $\hrc$ given by
$\hrc(a)=\ip{a z}z$, for all $a\in \arc$.  The Hopf
$*$-algebra associated to $(\arc,\ddrc)$ is $\aarc=\pi(\A)$, with 
co-unit $\e_{{\rm rc}}$ and co-inverse $\kappa_{{\rm rc}}$.
The map $\pi$ is a morphism of $\AA$ onto $(\arc,\ddrc)$ and its
kernel is equal to the left kernel of~$h$, so that $\pi$
induces a faithful representation of $\ar$ on~$H$. This
representation is an isomorphism of the compact quantum groups
$(\ar,\DD_{{\rm r}})$ and~$(\arc,\ddrc)$.
\end{thm}

\demo To prove that $\ddrc:\aarc\rightarrow\pi (\A )\otimes\pi
(\A )\subset B(H\otimes H)$ has an extension
$\ddrc:\arc\rightarrow B(H\otimes H)$, we construct a
unitary~$W$ on~${H\otimes H}$. First, define the linear map
${W:\A\otimes\A\subset H\otimes H\rightarrow
\A\otimes\A\subset H\otimes H}$ by setting $W(a\otimes
b)=\Delta (b)(a\otimes 1)$, for all $a,b\in\A$. We claim that
$W$ is isometric. To see this, let $c=\sum_i a_i\otimes
b_i\in\A\otimes\A$ and  $\Delta (b_i )=\sum_k a_i^k\otimes
b_i^k$ for finitely many elements $a_i^k ,b_i^k\in\A$. Then
\begin{eqnarray*}
\ \ &&W(c)^* W(c)= \sum_{ij} (\Delta (b_i )(a_i\otimes 1))^*
\Delta (b_j )(a_j\otimes 1)\\&&=\sum_{ijkl}(a_i^*\otimes 1
)((a_i^k)^*a_j^l\otimes (b_i^k )^* b_j^l )(a_j\otimes 1)
=\sum_{ijkl}a_i^*(a_i^k)^*a_j^l a_j\otimes (b_i^k )^*
b_j^l,\end{eqnarray*} and therefore \[\|W(c)\|_2^2
=(W(c)|W(c))=(h\otimes h)(W(c)^* W(c))\]
\[=\sum_{ijkl}h(a_i^*(a_i^k)^*a_j^l a_j )h((b_i^k )^* b_j^l )
=\sum_{ij}h(a_i^*[\sum_{kl}(a_i^k)^*a_j^l h((b_i^k )^* b_j^l
)]a_j )\] \[=\sum_{ij}h(a_i^*[(id\otimes h)\Delta (b_i^* b_j
)]a_j ) =\sum_{ij}h(a_i^* h(b_i^* b_j )1 a_j )
=\sum_{ij}h(a_i^* a_j )h(b_i^* b_j )\] \[=(h\otimes
h)(\sum_{ij}a_i^* a_j\otimes b_i^* b_j ) =(h\otimes h)(c^* c
)=(c|c)=\|c\|_2^2 .\] Hence $W$ is isometric, as claimed.
Since ${\A\otimes \A}$ is equal to the linear span
of~${\DD\A(\A\otimes 1)}$, we have ${W(\A\otimes
\A)}={\A\otimes \A}$. It follows that $W$ extends from the
dense subspace ${\A\otimes \A}$ to a unitary on~${H\otimes
H}$. We shall denote this extension also by~$W$.

We claim that, for all $a\in\A$, \begin{equation} \label{eqn:
W-equation} \ddrc (\pi (a))=W (\pi (1)\otimes\pi
(a))W^*;\end{equation} equivalently, $\ddrc (\pi (a))W=W(\pi
(1)\otimes\pi (a))$. These operators are equal if they act
identically on elementary tensors of the dense subspace
$\A\otimes\A$ of $H\otimes H$. Thus, let $b,c\in\A$ and
observe that \[\ddrc (\pi (a))W (b\otimes c)=\ddrc(\pi
(a))\Delta (c)(b\otimes 1) =((\pi\otimes\pi )\Delta (a))\Delta
(c)(b\otimes 1)\]
\[=\Delta (a)\Delta (c)(b\otimes 1)=\Delta (a c)(b\otimes 1)
=W(b\otimes a c)=W(\pi (1)\otimes\pi (a))(b\otimes c).\] Thus,
Equation~(\ref{eqn: W-equation}) holds and it follows that
$\ddrc:\aarc\rightarrow\arc\otimes\arc\subseteq B(H\otimes H)$
is norm decreasing. Consequently, it admits a $*$-homomorphism
extension ${\ddrc:\arc\rightarrow \arc\otimes \arc}$. That
$\ddrc$ is a co-multiplication on $\arc$ is an obvious
consequence of its restriction to $\aarc$ being one, and
density of $\aarc$ in~$\arc$. It follows directly, from the
fact that  the linear spans of ${(\A\otimes 1)\DD \A}$ and
${(1\otimes \A)\DD \A}$ are each equal to ${\A\otimes \A}$,
that $(\arc,\ddrc)$ is a compact quantum group. That $\pi$ is
a morphism of compact quantum groups is obvious.

Since $h=0$ on $\ker(\pi)$, it induces a unique state $\hrc$
on $\arc$ for which ${\hrc \circ \pi=h}$. Therefore, $\hrc(a)=\ip{a
z}z$, for all $a\in \arc$. It is easily verified that $\hrc$
is the  Haar state of $(\arc,\ddrc)$.  Suppose now $a\in \arc$
and ${\hrc (a^* a)=0}$. Since $N_{\hrc}$ is  a two-sided ideal
in $\arc$, $a b\in N_{\hrc}$, for all $b\in \arc$ and
therefore $\hrc (b^* a^* a b)=0$. Hence, $\ip{a b z}{a b
z}=0$ for all $b \in \arc$, which shows that $a=0$. Thus, $\hrc$ is
faithful. It
clearly follows that the left kernel of~$h$ is equal to the
kernel of~$\pi$.
Hence, the representation of $\ar$ on $H$ induced by $\pi$ is
faithful and then, by construction, an isomorphism of
$(\ar,\DD_{{\rm r}})$ onto $(\arc,\ddrc)$.

Finally, it is clear that $\pi(\A)$ is a
dense Hopf $*$-subalgebra of $(\aarc,\ddrc)$ with co-unit
$\e_{{\rm rc}}$ and co-inverse $\kappa_{{\rm rc}}$, and therefore it
is the Hopf $*$-algebra associated to~$(\aarc,\ddrc)$, by uniqueness .~\qed

We turn now to an application of some of our results to the
the prototypical example of a compact quantum group, the
quantization of $SU(2)$ constructed by Woronowicz. We shall
show that it is co-amenable, from which we shall obtain the
known, and non-trivial, result that its Haar integral is
faithful. It also follows from Banica's more general
result~\cite[Corollary 6.2]{Ba} which uses the theory of $\cal
R ^+$-deformations. Our quite elementary proof is totally
different.

Let $q$ be a real number for which $0<\m q< 1$. Let
$\AA=SU_q(2)$, and let $\a$ and $\g$ be the canonical
generators of $A$, satisfying the conditions of Table~0 of
\cite{Wo1}. Let $k\in\bZ$ and $m,n\in \bN$. Set
$a_{kmn}=\a^{(k)}\g^m\g^{*n}$, where $\a^{(k)}=\a^k$, if $k\ge
0$ and $\a^{(k)}=(\a^{-k})^*$, if $k<0$. Recall that these
elements $a_{kmn}$ form a linear basis for the Hopf *-algebra
$\sa$ associated to $\AA$ and that $h(a_{kmn})=0$, if $k\ne 0$
or if $m\ne n$~\cite[Equation~A1.8]{Wo1}.

Take $U$ to be the fundamental irreducible co-representation
of~$SU_q(2)$ given by
\[U=\left(\begin{array}{cc}
 \a &-q\g^*\\ \g&\a^*\end{array}\right).\]

Before stating the following theorem, we make an elementary
observation: If $V$ is the forward unilateral shift on a
Hilbert space $H$ with orthonormal basis $(e_n)_{n\in \bN}$,
so that $Ve_n=e_{n+1}$, then there exists a state $\t$ on
$B(H)$ such that $\t(V)=1$ and $\t(K)=0$, for all compact
operators $K\in B(H)$. To see this, one observes that the
image of $V$ in the Calkin algebra~$C$ of $H$ is a unitary
containing 1 in its spectrum and therefore there exists a
state on $C$ whose value at this unitary is equal to~1. The
required state on $B(H)$ is then the composition of the state
on $C$ and the quotient map from $B(H)$ to~$C$.

\begin{thm} The compact quantum group $SU_q(2)$ is co-amenable.
\end{thm}

\demo As before, let $(A,\DD)=SU_q(2)$ and let $\a$ and $\g$
be the canonical generators of $A$. Set
${c_n=(1-q^{2n})^{1/2}}$, for $n\in\bN$.  Recall from
Appendix~A.1 of~\cite{Wo2} that~$A$ admits a representation
$\f$ on a Hilbert space $H$ with an orthonormal basis
$(e_{n,k})$, where $n\in\bN$ and $k\in \bZ$, such that
\[ \f(\a)e_{nk}=c_ne_{n-1,k}\quad {\rm and }\quad \f(\g)e_{nk}=q^n
e_{n,k+1}\] and that
\[h(a)=(1-q^2)\sum_{n=0}^\infty \ip{\f(a)e_{n0}}{e_{n0}}.\]
It follows immediately that $h(a^*a)=0$ if, and only if,
$\f(a)e_{n0}=0$, for all $n\in \bN$. Using the equations
$\f(\g^m)e_{n0}=q^{nm}e_{nm}$ and
$\f(\g^{*m})e_{n0}=q^{nm}e_{n,-m}$, for $m>0$ and the fact
that $a\g^m$ and $a\g^{*m}$ belong to $N_h$, if $a$ does, we
get that $h(a^*a)=0$ if, and only if, $\f(a)=0$. Hence, we get
an induced faithful representation $\p$ of $\ar$ on $H$ given
by $\p\ta(a)=\f(a)$.

Now, for $k\in \bZ$, let $H_k$ be the Hilbert subspace of $H$
with orthonormal basis $(e_{nk})_{n\in \bN}$. Obviously,
$H=\oplus_k H_k$, and $T=\f(\a)$ reduces each space~$H_k$, so
that $T=\oplus_k T_k$, where $T_k$ is the restriction of $T$
to~$H_k$. We have $T_ke_{nk}=c_ne_{n-1,k}$, so that
$T_k=U^*_kD_k$, where $U_k$ is the forward unilateral shift on
the basis $(e_{nk})_n$ of $H_k$ and $D_k$ is the diagonal
norm-bounded linear operator on~$H_k$ defined by setting
$D_k(e_{nk})=c_ne_{nk}$. Since $\lim c_n=1$, it is clear that
$D_k=1+L_k$, where $L_k$ is a compact operator on~$H_k$.
Hence,  $T_k=U^*_k+U^*_kL_k$. By the remarks preceding this
theorem, there exists a state $\t_\k\in B(H_k)$ such that
$\t_k(T_k)=1$. For $k\in \bZ$, chose positive numbers $t_k$
such that ${\sum_{k\in \bZ} t_k=1}$. Now define a state $\t$
on the C*-algebra $\oplus_k B(H_k)$ containing $T$ by setting
$\t(S)=\sum_{k\in \bZ} t_k \t_k(S_k)$, if $S=(S_k)_k\in
\oplus_k B(H_k)$. Clearly, $\t(T)=1$. Now let $\t'$ be the
state $\t\p$ on $\ar$. Then $\t'(\ta(\a))=\t(T)=1$ and
therefore $\t'\ta(\Re \chi_U)=\t'\ta(\a)+(\t'\ta(\a))^-=2$.
Hence, $\AA$ is co-amenable, by Condition~(3) of
Corollary~\ref{thm: co-amenability corollary}.~\qed

\begin{cor}[G. Nagy] \label{thm; Nagy's theorem}
The Haar integral $h$ of $SU_q(2)$ is
faithful.
\end{cor}

\demo This is a consequence of the preceding theorem and
Theorem~\ref{thm: Haar faithfulness}.~\qed

There is an alternative way of proving $SU_q(2)=(A,\DD)$ is
co-amenable, using the fact that $A$ is of Type~I, as a
C*-algebra~\cite[Theorem~A2.3]{Wo1}. Since $\ar$ is unital, it
admits a maximal ideal $I$. Since $\ar/I$ is a Type~I simple
unital C*-algebra, it is isomorphic to $M_N(\bC)$, for some
positive integer~$N$. Thus, we have a surjective
$*$-homomorphism $\pi$ from $\ar$ onto~$M_N(\bC)$. The
existence of a faithful, tracial state on $M_N(\bC)$, together
with the commutation relations of~\cite[Table~0]{Wo1} for the
canonical generators $\a$ and $\g$, forces the image $\pi(\g)$
of $\g$ in $M_N(\bC)$ to be equal to zero and $\pi(\a)$ to be
a unitary. Since $\pi(\a)$ and $\pi(\g)$ generate $M_N(\bC)$,
this implies that $M_N(\bC)$ is commutative. Hence, $N=1$ and
$M_N(\bC)=\bC$. Thus, $\ar$ admits a $*$-homomorphism
onto~$\bC$ and it now follows from~Corollary~\ref{thm:
character corollary} that $SU_q(2)$ is co-amenable.

\section{The Universal Quantum Group}

In this section we first give a detailed account on the construction
of the universal compact quantum
group associated to an arbitrary compact quantum group.
One way to construct such an object relies on Baaj and Skandalis'
theory of regular multiplicative unitaries \cite{BS}.
A general construction for locally compact quantum groups has
recently been given by J.~Kustermans~\cite{Ku}. However, our
approach, which is briefly sketched by Woronowicz in
\cite{Wo3} for compact  matrix pseudogroups, is much less
technical and is therefore included. 
The reduced quantum group has the advantage that the Haar integral
is always faithful, whereas its co-unit need not be
norm-bounded. For the universal quantum group the situation is the opposite; 
its co-unit is always norm-bounded, whereas its Haar integral
need not be faithful.

\vspace{2ex}

Let $\AA$ be a compact quantum group. Define 
$\nuu{\cdot}$ on $\A$ by
\[\nuu a = \sup_\pi \n{\pi(a)},\]
where the variable $\pi$ runs over all unital
$*$-homomorphisms $\pi$ from $\A$ into $B(H_\pi)$, for a
Hilbert space~$H_\pi$ (the {\em unital $*$-representations} of
$\A$).

\begin{lem}
The function $\|\cdot\|_u :\A\rightarrow [0,\infty ]$ is a
$C^*$-norm on $\A$ which majorises any other $C^*$-norm on
$\A$.
\end{lem}

\demo We first need to show that $\nuu a$ is finite, for all
$a\in\A$. Let $(U^\a)_\a$ be a complete set of inequivalent,
irreducible unitary co-representations of~$\AA$; then the
matrix elements $U^\a_{ij}$ linearly span $\A$. Clearly, it
suffices to show that $\nuu{U^\a_{ij}}<\infty$, for all $\a$
and $i,j$. Suppose then ${\pi :\A\rightarrow B(H)}$ is a
unital $*$-representation of $\A$ on some Hilbert space~$H$.
Since $\sum_k (U^{\alpha}_{ki})^* U^{\alpha}_{ki} =1$, we have
$(U^{\alpha}_{ji})^* U^{\alpha}_{ji}=1-\sum_{k\neq j}
(U^{\alpha}_{ki})^* U^{\alpha}_{ki}$, and therefore
\[0\leq (\pi (U^{\alpha}_{ji}))^*  \pi (U^{\alpha}_{ji}) =\pi
(1)-\sum_{k\neq j} (\pi (U^{\alpha}_{ki}))^* \pi
(U^{\alpha}_{ki} )\leq \pi (1).\] Hence, $\|\pi
(U^{\alpha}_{ji} )\|^2 =\|(\pi (U^{\alpha}_{ji}))^*  \pi
(U^{\alpha}_{ji})\|\leq \|\pi (1)\|=1$. It follows that
$\nuu{U^\a_{ij}}$ is finite. (Note that although $0\leq\|a^*
a\|I-a^* a$, for any $a\in \A$, we cannot conclude $\|a^*
a\|I-a^* a =b^*b$, for some element $b$ belonging to $\A$.
This is why the preceding argument had to be more careful than
one might first expect  and had to use the rather strong
property that $\A$ is the linear span of the matrix
elements~$U^\a_{ij}$.)

It is clear now that $\nuu{\cdot}$ is a C*-seminorm on~$\A$
and since $A$ admits a faithful unital $*$-representation,
$\nuu{\cdot}$ is, in fact, a C*-norm. That $\nuu{\cdot}$
majorises any other C*-norm on~$\A$ is clear from its
definition.~\qed

We define $\au$ to be the $C^*$-algebra completion of $\A$
with respect to the $C^*$-norm $\nuu{\cdot}$. As usual, we identify 
$\A$ with its canonical copy inside $\au$. The C*-algebra $\au$ has
the universal property that if ${\pi\colon \A\to B}$ is a
unital $*$-homomorphism from $\A$ to a unital C*-algebra~$B$,
it extends uniquely to a $*$-homomorphism from $\au$ to~$B$, since
$\pi$ is easily seen to be norm-decreasing on $\A$ equipped with its
universal norm. 

In particular, the $*$-homomorphism ${\DD\colon \A\to
\A\otimes \A\subseteq \au\otimes \au}$ extends to a
$*$-homomorphism ${\DD\colon \au\to \au\otimes \au}$. It is
easily verified $\DD$ is a co-multiplication on~$\au$. Since
the linear spans of the sets ${(\A\otimes 1)\DD\A}$ and
${(1\otimes \A)\DD\A}$ are each equal to~${\A\otimes \A}$, it
follows immediately that $\AAu$ is a compact quantum group. We
call it the {\em universal} compact quantum group associated
to~$\AA$.

\vspace{1ex}
Since $\A$ is, by construction, a dense Hopf $*$-subalgebra of 
$\AAu$, it is the Hopf $*$-algebra associated to $\AAu$, by uniqueness.

Note also that the co-unit $\e$ of $\A$, being a $*$-homomorphism from
$\A$ to $\bC$, extends to a $*$-homomorphism $\eu$ from
$\au$ to $\bC$.  By density of $\A$ in $\au$, the equalities
${(\e\otimes \id)\DD(a)}={(\id\otimes \e)\DD(a)}=a$, which
hold for all $a\in \A$, extend to ${(\eu\otimes
\id)\DD(a)}= {(\id\otimes \eu)\DD(a)}=a$, for all $a\in \au$.
Hence, $\eu$ must be the unique extension to $\au$ of the
co-unit of~$\AAu$. The important point we wish to emphasize
here is that $\AAu$ has thus a norm-bounded co-unit.

Clearly, by the universal property of $\AAu$, there is a
$*$-homomorphism $\p$ from $\au$ onto $A$ extending the identity
$*$-isomorphism from $\A$ to itself. Also, $\DD\p={(\p\otimes
\p)\DD}$. We call $\p$ the {\em canonical map} from $\au$ onto
$A$. Likewise, if $\ta$ is the canonical map from $A$ onto
$\ar$, we call the composition $\ta\p$ the {\em canonical map}
from $\au$ onto $\ar$.

 Clearly, $h\p$ is the Haar integral~$\hu$ of $\AAu$; hence,
$\hu=\hr\ta\p$. Since $\hr$ is faithful, it follows that
$N_{\hu}=\ker(\ta\p)$. From this it is immediate that the
reduced compact quantum group of $\AAu$ is (isomorphic to)
$(\ar,\DD_{{\rm r}})$ and that $\p\ta$ is the canonical map
from $\AAu$ onto
$(\ar,\DD_{{\rm r}})$. Therefore, $\AAu$ is co-amenable if,
and only if,
$\AA$ is co-amenable.

We summarise the preceding discussion in the following
theorem.

\begin{thm} \label{thm: universality theorem}
Let $\AA$ be a compact quantum group. Then $\A$ is the Hopf 
$*$-algebra associated to the universal compact quantum group  
$\AAu$. The  co-unit of $\AAu$ is norm-bounded. 
Finally, the reduced compact quantum group of
$\AAu$ is (isomorphic to) $(\ar,\DD_{{\rm r}})$, so that  $\AAu$ is
co-amenable if, and only if, $\AA$ is.
\end{thm}

It is quite obvious that the universal compact quantum group
${(C^*(\Gamma),\DD)}$ associated to a discrete group $\Gamma$
is its own universal compact quantum group; that is, if
$\AA={(C^*(\Gamma),\DD)}$, then $\AAu=\AA$. Moreover, if
$\AA={(C^*_{{\rm r}}(\Gamma),\DD_{{\rm r}})}$, then
$\AAu={(C^*(\Gamma),\DD)}$. This is the
 motivating example for the general definition of
the universal compact quantum group.

Suppose now $\AA$ is an arbitrary compact quantum group. It is
easy to see that if $(B,\Phi)$ is a compact quantum group
whose associated Hopf $*$~-algebra $(\B,\Phi)$ is isomorphic
to $(\A,\DD)$, then $(B_{u},\Phi)$ is isomorphic to $\AAu$.
In particular, the universal compact quantum group associated
to $(\ar,\DD_{{\rm r}})$, or to $\AAu$, is isomorphic to $\AAu$.

We call a compact quantum group $\AA$ {\em universal} if $\AA=\AAu$, 
\\ i.\ e.\ if  the  canonical map $\p$ from $\au$ onto $A$ is 
injective. Equivalently, $\AA$ is universal if, and only if, the given 
norm on $\A$ is its greatest C*-norm. We will show in 
Corollary~\ref{pos} that any co-amenable compact quantum group is universal.

We prove now a striking automatic continuity result for
positive linear functionals on the Hopf $*$-algebra of a
universal compact quantum group. Recall that a linear 
functional $\t$ on a $*$-algebra $B$ is called {\em
positive} if $\t(b^*b)\ge 0$, for all $b\in B$.

\begin{thm}
\label{thm: functional boundedness theorem} Suppose that
$(A,\Delta)$ is a universal compact quantum group. Then
every positive linear functional $\t$ on $\A$ is norm-bounded.
\end{thm}

\demo  We form the GNS representation of $\A$ with
respect to $\t$: Since the map ${(a,b)\mapsto
\t(b^*a)}$ is sesquilinear, the inequality
$\m{\t(b^*a)}^2\le \t(b^*b)\t(a^*a)$ implies that the left
kernel $N_\t$ of $\t$ is a left ideal of $\A$. Hence, the
quotient space $\A/N_\t$  is a inner product space with inner
product given by $\ip{a+N_\t}{b+N_\t}=\t (b^* a)$, where
$a,b\in\A$. Denote the Hilbert space completion by $H$ and its
norm by $\|\cdot\|_2$. Define the operator
$M_a:\A/N_\t\to\A/N_\t$ by setting ${M_a (b+N_\t)}={a
b+N_\t}$, for all $a,b\in\A$.

We shall show now that $M_a$ is norm-bounded, for all $a\in
\A$. Since the map, ${a\mapsto M_a}$, is linear, it suffices
to show boundedness for $a=U^\a_{ij}$, where $(U^\a)_\a$ is a
complete set of inequivalent, irreducible unitary
representations of~$\AA$ and $U^\a_{ij}$ are the matrix
elements of $U^\a$. We have, for all $b\in\A$,   \[{b^* b
-b^*(U^{\alpha}_{ji})^* U^{\alpha}_{ji}b}={b^* (\sum_{k\neq
j}(U^{\alpha}_{ki})^* U^{\alpha}_{ki})b} ={\sum_{k\neq
j}(U^{\alpha}_{ki}b)^* (U^{\alpha}_{ki}b)\geq 0}.\] Hence,
$\|U^{\alpha}_{ji}b+N_\t\|_2^2 =\t(b^*(U^{\alpha}_{ji})^*
U^{\alpha}_{ji}b)\leq \t(b^* b)=\|b+N_\t\|_2^2$, so that
${\n{M_a}\le 1}$. (This kind of argument was used tacitly in
the proof of the generalized Tannaka-Krein theorem
in~\cite{Wo3}.)

Hence, for all $a\in\A$, we may extend $M_a$ to a
norm-bounded operator $\pi (a)$ on $H$. The corresponding map,
${\pi:\A\rightarrow B(H)}$, ${a\mapsto \pi(a)}$,  is obviously
a unital $*$-representation of $\A$. By the universal property
of $\au$, this map extends to a $*$-homomorphism
${\pi\colon \au\to B(H)}$. 
Since, for all $a\in\A$, $\t
(a)=\ip{\pi (a)x}x$, where $x={1+N_\t}$, we have \[ \m{\t
(a)}=\m{\ip{\pi(a)x}x}\leq \|\pi (a)\|\|x \|_2^2=\|\pi (a)\|\t
(1^* 1)\leq \|a\|_u\t(1). \] Hence, $\t$ is norm-bounded with
respect to the universal C*-norm on~$\A$. Since $\AA$ is assumed to 
be universal, this norm is equal to the given norm on~$\A$.~\qed

\vspace{2ex}
When $A$ is a unital C*-algebra, one may consider the
C*-algebra invariant consisting of all non-zero
$*$-homomorphisms from $A$ to $\bC$, i.\ e.\ of all unital multiplicative
linear functionals on $A$. This (possibly empty) set
is clearly compact in the relative weak* topology inherited
from $A^{*}$. Of course, when $A$ is commutative, it is precisely the
Gelfand spectrum of $A$. For some other classes of (non-simple)
$C^{*}$-algebras, this generally rather poor invariant is of
some interest. For example, when $A$ is the universal compact group 
associated to a discrete group $\Gamma$, it is easily
identified with the dual group of the abelianized group of
$\Gamma$ (see \cite{Wat}) and therefore it is
computable in many cases. We will show below that this
invariant is a compact group for any universal compact quantum
group. 

\vspace{1ex}
We need a lemma which may be known to specialists, but for
which we could not find a suitable reference in the
literature.

\begin{lem} \label{thm: group} If $\AA$ is a compact quantum group,
 the unital multiplicative linear functionals on~$\A$ form
a group under the multiplication, ${(\t,\s)\mapsto \t*\s}$,
where $\t*\s={(\t\otimes\s)\DD}$. The unit is $\e$ and the
inverse of the element $\t$ is $\t\k$. Moreover, the
$*$-homomorphisms from $\A$ onto $\bC$ form a subgroup (which
may be proper).

\end{lem}

\demo That the operation is closed and associative and the
co-unit is a unit for this operation is well known. We prove
first that the inverse of the element $\t$ is $\t\kappa$. To
see $\t*(\t\k)=\e$, let $a\in \A$, and observe that
$(\t*(\t\kappa))(a)={(\t\otimes \t\k)\DD(a)} ={\t(
m(\id\otimes \k)\DD(a))}=\t(\e(a)1)=\e(a)$. Here ${m\colon
\A\otimes \A\to \A}$ is the linearization  of the
multiplication ${\A\times \A\to \A}$. We used the fact that
${\t\otimes \t\k}={\t m(\id\otimes \k)}$ which is a
consequence of the multiplicative property enjoyed by~$\t$.
That $(\t\k)*\t=\e$ is similarly proved. Now if
${\t\colon\A\to \bC}$ is a $*$-homomorphism, then $\t\k$ is
also. We prove this indirectly. The map $\hat
\t={(\id\otimes\t)\DD\colon\A\to \A}$ is a $*$-homomorphism,
since $\t$ is one. Moreover,
${\hat\t((\t\k)\hat{\,}(a))}={(\t*\t\k)\hat{\,}(a)}=\hat\e(a)=a$
and likewise
${(\t\k)\hat{\,}(\hat{\t}(a))}={((\t\k)*\t)\hat{\,}(a)}=\hat\e(a)=a$.
Hence, $(\t\k)\hat{\,}$ is the inverse of $\hat\t$ and it is
therefore also a $*$-homomorphism. Finally, since
${\t\k=\e\circ(\t\k)\hat{\,}}$ is a composition of
$*$-homomorphisms, it is one also. Hence,  the
$*$-homomorphisms from $\A$ onto $\bC$ form a subgroup, which
may be proper since  multiplicative linear functionals on a *-algebra
do not necessarily preserves adjoints.~\qed

\begin{thm} \label{thm: character group theorem}
If $\AA$ is a universal compact quantum group,
then the set $G$ of unital multiplicative linear functionals
on~$A$ forms a compact topological group under the relative
weak* topology and the multiplication, ${(\t,\s)\mapsto
\t*\s}$, where $\t*\s={(\t\otimes\s)\DD}$.
\end{thm}

\demo As before, closure and associativity of the
multiplication operation is well known and since the co-unit
of $\AA$ is norm-bounded, its extension to $A$ exists and
provides a unit for~$G$. If $\t$ is a unital multiplicative
linear functional on $A$, it is necessarily a
$*$-homomorphism. Hence, if $\t$ is its restriction to $\A$,
the functional $\t\k$ is also a $*$-homomorphism, by the
preceding lemma. By universality of~$\AA$, $\t\k$ admits an
extension to a $*$-homomorphism, $\s$ say, on~$A$. Since
${(\t\otimes \s)\DD(a)}={(\s\otimes \t)\DD(a)}=\e(a)$, for all
$a\in \A$, the same equalities hold for all $a\in A$, by
continuity. Thus, $\t*\s=\s*\t=\e$. It is straightforward to
check that $G$ is a weak* closed subset of the unit ball
of~$A^*$ and therefore, by  the Banach--Alaoglu theorem, $G$
is weak* compact. It is also easily checked that the
multiplication operation is weak* continuous, as is the
inversion operation ${\t\mapsto\t^{-1}}$. This proves the
theorem.~\qed

As an example,
let $\AA$ be the compact quantum group $SU_q(2)$, where $q\in
\bR$ and $0<\m q< 1$. Being co-amenable, $\AA$ is universal.
Let $\a$ and $\g$ be the canonical generators of~$A$. If $\t$ belongs
to the group $G$ of multiplicative linear functional on $A$,
then the equations $\a\a^*+\g\g^*=1=\a^*\a+q^2\g^*\g$ imply
that $\t(\g)=0$ and $\t(\a)$ belongs to the unit circle
group~$\bT$. Conversely, given $\l\in \bT$, the universal
property enjoyed by $A$ implies that there exists
a---necessarily unique---element~$\t$ of $G$ for which
$\t(\a)=\l$ (and $\t(\g)=0$). Since
$\DD\a={\a\otimes\a-q\g^*\otimes \g}$, we have
$(\t*\s)(\a)=\t(\a)\s(\a)$, for all $\t,\s\in G$. Hence the
map, ${\t\mapsto \t(\a)}$,  is a group isomorphism from $G$
onto~$\bT$. It is trivially continuous, so that it is also a
homeomorphism (since the spaces are compact and Hausdorff).
Thus, $G=\bT$, as topological groups.

Lemma~\ref{thm: group} can be used to give an alternative
proof of Corollary~\ref{thm: character corollary}: Let $\AA$
be a compact quantum group and suppose given a
$*$-homomorphism ${\t\colon A\to \bC}$. Of course, its
restriction to $\A$ is therefore a $*$-homomorphism, from
which it follows that $\t\k$ is one also. Hence, by
\cite[Lemma~10.2]{MT}, $(\t\k)\hat{\,}$ is an isometry (we
are retaining the notation used in the proof of
Lemma~\ref{thm: group}). Since $\e=\t\circ(\t\k)\hat{\,}$ is
the composition of two norm-bounded maps, it is norm-bounded
and therefore $\AA$ is co-amenable.

\vspace{2ex}
We now come to one of the main results of the theory. It especially
confirms that the Haar integral of a co-amenable compact quantum group
is faithful. The equivalence between (1) and (2) shows that our definition
of co-amenability agrees with the one considered by Banica \cite{Ba,Ba2}.

\begin{thm}
\label{am} The following are equivalent conditions for a
compact quantum group $(A,\Delta)$:

(1) $\AA$ is co-amenable;

(2) The canonical map from $\au$ to $\ar$ is a
$*$-isomorphism;

(3) The canonical maps from $\au$ onto $A$ and $A$ onto
$\ar$ are $*$-isomorphisms;

(4) The Haar integral $h_u$ of $\AAu$ is faithful.
\end{thm}

\demo If Condition~(1) holds, then $\AAu$ is co-amenable, by
Theorem~\ref{thm: universality theorem} and therefore $\hu$ is
faithful, by Theorem~\ref{thm: Haar faithfulness}. Thus,
$(1)\implies (4)$. Since $\hu=h\p=\hr\ta\p$, it is clear that
Condition~(4) implies~(2). The equivalence of Conditions~(2)
and~(3) is trivial. Suppose now that~(2) holds and let $\eu$
be the extension of the co-unit of $\AAu$ to $\au$
Then $\eu(\ta\p)^{-1}$ is a non-zero $*$-homomorphism on $\ar$ and
therefore, by Corollary~\ref{thm: character corollary}, $\AA$
is co-amenable. Thus, $(2)\implies (1)$. This proves the
theorem.~\qed

The following is now immediate from the theorem, from 
Theorem~\ref{thm: functional boundedness theorem}
and from Theorem~\ref{thm: character group theorem}.

\begin{cor} \label{pos}
Let $\AA$ be a co-amenable compact quantum group. Then $\AA$ is 
universal. Especially, every unital
$*$-homomorphism from $\A$ to a unital C*-algebra  is
necessarily norm-decreasing. Further, every positive linear 
functional on $\A$ is norm-bounded. Finally, 
the unital multiplicative linear functionals on~$A$ form a compact
group 
\end{cor}

Note that co-amenability imposes a norm-boundedness condition on just a
single postive linear functional (the co-unit of the reduced
quantum group). However, the corollary shows it
implies a much stronger norm-boundedness result.

The equivalence between (1) and (3) in Theorem~\ref{am} may be 
rephrased as saying that a compact quantum group $\AA$  is co-amenable 
if, and only if, it is both universal and reduced. 
Note in this connection that $C^*({\bf F}_2) \otimes C^*_{r}({\bf F}_2)$ is
an example of a compact quantum group which is neither universal nor
reduced, since, obviously, its Haar integral is not faithful and its
co-unit is not norm bounded.

\vspace {2ex}
If $\AA$ is an arbitrary compact quantum group, we know
that $\nuu{\cdot}$ is the greatest C*-norm on the associated
Hopf $*$-algebra~$\A$. We define a C*-seminorm on~$\A$ by
setting $\nr{a}=\n{\ta(a)}$, for all $a\in\A$. This is, in
fact, a C*-norm, since $\ta$ is injective on~$\A$. Therefore
we can regard not only $\AAu$ and $\AA$ as compact quantum
group completions of~$\A$, but $(\ar,\DD_{\rm r})$ also.
When we say that a compact quantum group
${(A_{\rom c},\DD_{\rom c})}$ is a {\em compact quantum group
completion} of $\A$, we mean not only that $\A$ is a
dense unital $*$-subalgebra of the C*-algebra, but also that
the co-multiplication $\DD_{\rom c}$ extends the
co-multiplication $\DD$ of $\A$. We shall call a C*-norm
$\nc{\cdot}$ on~$\sa$ {\em regular} if it is the restriction
to $\A$ of the norm of a compact quantum group completion
${(A_{\rom c},\DD_{\rom c})}$ of~$\A$. Thus, the given C*-norm on~$\A$
and the norms $\nuu{\cdot}$ and $\nr{\cdot}$ are regular.

\vspace {1ex}
We show now that $\nr{\cdot}$ is the least regular C*-norm
on~$\A$.

\begin{thm} \label{thm: norm inequalities}
Let $\AA$ be a compact quantum group and $\nc{\cdot}$ be a
regular C*-norm on~$\A$. Then $\nr a\le \nc a\le \nuu a$, for
all $a\in \A$. If ${(A_{\rom c},\DD_{\rom c})}$ is the compact quantum
group completion of~$\A$ with respect to~$\nc{\cdot}$, then
there exist unique $*$-homomorphisms ${\p_{\rom c}\colon
\au\to A_{\rom c}}$  and ${\ta_{\rom c}\colon A_{\rom c}\to
\ar}$ extending, in each case, the identity automorphism on
$\A$. Both maps are quantum group morphisms.
\end{thm}

\demo Given the maps $\p_{\rom c}$ and $\ta_{\rom c}$ exist,
it follows trivially from density of~$\A$ in~$\au$ and
$A_{\rom c}$, respectively, that they are unique and are
quantum group morphisms. The
norm inequality $\nc{\cdot}\le \nuu{\cdot}$ is already known and the existence
of the map $\p_{\rom c}$ is obvious. If we show $\nr{\cdot}\le
\nc{\cdot}$, the existence of $\ta_{\rom c}$ follows
trivially. We turn now to showing this inequality. Before
proceeding, let us first observe that $\m{h(a)}\le \nc a$, for
all $a\in \A$. Let $h_{\rom c}$ denote the Haar integral of  
${(A_{\rom c},\DD_{\rom c})}$. When we regard $\A$ as a  Hopf $*$-subalgebra
of ${(A_{\rom c},\DD_{\rom c})}$ and of $\AA$, as we do here, we have
$h_{\rom c}(a) = h(a)$ for all $a \in \A$, by uniqueness of    
Haar integrals. 
Consequently, $\m{h(a)}=\m{h_{\rom c}(a)}\le \nc a$, as
claimed.

Again suppose $a\in \A$. Since the Haar integral~$\hr$
of~$(\ar,\DD_{\rm r})$ is a faithful state of~$\ar$, it follows from
\cite[Theorem~10.1]{MT} that
\[\nr{a^*a}=\n{\ta(a)^*\ta(a)}=
{\lim[\hr((\ta(a)^*\ta(a))^n)]^{1/n}}.\] Using the fact that
$h=\hr\ta$, we get $\nr{a^*a}={\lim[h(a^*a)^n)]^{1/n}}$.  By our
observations in the preceding paragraph, ${h((a^*a)^n)\le
\nc{(a^*a)^n}}$. Therefore, $\nr{a^*a}\le
{\lim\nc{(a^*a)^n}^{1/n}}=\nc{a^*a}$ and hence $\nr a\le \nc
a$, as required.~\qed

\begin{cor}
Let $\AA$ be a compact quantum group with associated Hopf
$*$-algebra $\A$. Then $\AA$ is co-amenable if, and
only if, $\A$ admits only one quantum group completion
(up to isomorphism).
\end{cor}

\demo This follows immediately from the theorem and the
observation that $\AA$ is co-amenable if, and only if,
$\nuu{\cdot}=\nr{\cdot}$, as norms on $\A$; this observation
is an immediate consequence of Theorem~\ref{am}.~\qed

The qualifying word {\em regular} may not be dropped in the
statement of Theorem~\ref{thm: norm inequalities}. This may be
seen as follows: Let $\Gamma$ denote a discrete group. Set
$\AA=(C^{*}(\Gamma), \DD)$, and recall that $\A$ is the group
algebra of $\Gamma$. Let $W$ be a unitary representation of
$\Gamma$ on a Hilbert space $H$ and denote by $\pi$ the
associated representation of $C^{*}(\Gamma)$ on $H$, so that
$\pi(C^{*}(\Gamma))=C^{*}(W)$,
where $C^{*}(W)$ denotes the $C^{*}$-algebra generated by all
$W_{x}$ ( $x \in \Gamma $). Then define a $C^{*}$-seminorm
$\np{\cdot}$ on $\A$ by setting $\np a = \n{\pi(a)}$. Assume
that $\np{\cdot}$ is a $C^{*}$-norm on $\A$; that is, $\pi$ is
faithful on $\A$. Then the completion of $\A$ with respect to
$\np{\cdot}$ may be identified with $C^{*}(W)$.

If we now assume that Theorem~\ref{thm: norm inequalities}
holds without the qualifying word {\em regular}, the regular
representation $L$ of $\Gamma$ is clearly weakly contained in
$W$; that is, there exists a $*$-homomorphism $\phi $ from
$C^{*}(W)$ onto $C^{*}(L)$ satisfying $\phi(W_{x})= L_{x}$,
for all $ x \in \Gamma$. If we also assume that $\Gamma$ is
amenable, then $\f$ is a $*$-isomorphism (since it clearly
admits an inverse in this case). Now set $\Gamma=\bZ$. Then
$C^*(L)=C(\bT)$ and $L_1$ has spectrum $\bT$. This forces
$W_1$ to have spectrum~$\bT$ also. To get a contradiction we
need now only show $W_1$ does not have to have spectrum~$\bT$.
To do this, choose a unitary $V$ on a Hilbert space with
infinite spectrum not equal to~$\bT$. This induces a
representation $W$ of $\bZ$ and the corresponding homomorphism
$\pi$ is injective on ${\bf C}(\bZ)$, since ${{\rm s p}}(V)$
is infinite (this implies all the powers ${1,V,V^2,\dots}$ are
linearly independent). Thus, this representation $W$ satisfies
the required conditions and the spectrum of $W_1=V$ is not
equal to~$\bT$.

An open question in this setting is whether  $\np{\cdot}$ is
necessarily regular whenever $L$ is weakly contained in $W$.
We doubt that the answer is positive.
It is worth mentioning here that
Woronowicz shows in \cite[Theorem 1.6]{Wo1} that if $\Gamma$
is finitely generated and $W$ is a faithful representation of
$\Gamma$ such that $W \otimes W$  is (strongly) contained in a
multiple of $W$, then $\np{\cdot}$ is regular.
However, the only known representations satisfying these
assumptions seem to be the universal and the regular ones, and the
external tensor product of such representations.

\section{Quantum semigroups and co-amenability}

In this short section we give a sufficient condition 
ensuring that a compact quantum semigroup is a compact quantum
group. Recall that a {\em compact quantum semigroup} is a pair
$\AA$ consisting of a unital C*-algebra~$A$ and a
co-multiplication ${\DD\colon A\to A\otimes A}$. Of course,
if, in addition, the linear spans of the spaces ${(A\otimes
1)\DD A}$ and ${(1\otimes A)\DD A}$ are each equal to
${A\otimes A}$, then $\AA$ is a compact quantum group. A {\em
Haar integral} on a compact quantum semigroup $\AA$ is defined
in the usual way as a state on $h$ on~$A$ for which we have
${(\id\otimes h)\DD(a)}={(h\otimes \id)\DD(a)}=h(a)1$, for all
$a\in A$. It is trivial to verify that at most one Haar
integral can exist. Not every compact quantum semigroup admits
a Haar integral, nor does  the existence of a Haar integral
imply that a compact quantum semigroup is a compact quantum
group~\cite{MT}.

A {\em bounded co-unit} for a compact quantum semigroup $\AA$
is defined as a unital $*$-homomorphism $\e$ from $A$ to $\bC$
such that, for all $a\in A$, ${(\e\otimes
\id)\DD(a)}={(\id\otimes \e)\DD(a)}=a$. The example given
in~\cite{MT} of a compact quantum semigroup having no Haar
integral has got a bounded co-unit. Thus, the existence of a
bounded co-unit does not ensure that a compact quantum
semigroup is a compact quantum group.

We shall need some notation for the following two results. If
$a,b\in A$, we write $a*(h b)$ for the element ${(h\otimes
\id)((b\otimes 1)\DD(a))}$ and $(h a)*b$ for the element
${(\id\otimes h)((1\otimes a)\DD(b))}$.

\begin{lem}
Let $\AA$ be a compact quantum semigroup admitting a Haar
integral~$h$. Then, for all $a,b\in A$, the element ${1\otimes
a*hb}$ belongs to the closed linear span of ${(A\otimes 1)\DD
A}$. Likewise, ${(ha)*b\otimes 1}$ belongs to the closed
linear span of ${(1\otimes A)\DD A}$.
\end{lem}

\demo If $F\colon A\otimes A\to A\otimes A$ is the flip
automorphism, then the {\em opposite} compact quantum
semigroup $(A,F\DD)$ also has the state $h$ as its Haar
integral and $\e$ as a bounded co-unit. It follows that if we
show that ${1\otimes a*hb}$ belongs to the closed linear span
of ${(A\otimes 1)\DD A}$, then we can deduce from this result
applied to $(A,F\DD)$ that  ${(ha)*b\otimes 1}$ belongs to the
closed linear span of ${(1\otimes A)\DD A}$. The demonstration
that ${1\otimes a*hb}$ belongs to the closed linear span of
${(A\otimes 1)\DD A}$ is given in the proof of Theorem~3.3 of
\cite{MT}. The strong hypotheses of Theorem~3.3 are not needed
for our result, which only needs the fact that $\AA$ is a
compact quantum semigroup admitting a Haar integral, as can be
verified by a careful reading of the proof in~\cite{MT}.~\qed

\begin{thm}
Let $\AA$ be a compact quantum semigroup admitting a faithful
Haar integral and a bounded co-unit. Then $\AA$ is a
co-amenable compact quantum group.
\end{thm}

\demo If we show that $\AA$ is a compact quantum group, its 
co-amenability follows from Theorem~\ref{thm: Haar faithfulness}.
By the preceding lemma, we need then only show that the
closed linear span $L$ of the elements $a*hb$, where $a,b \in
A$, and the closed linear span $R$ of the elements $(ha)*b$,
are both equal to~$A$. For, in this case, ${1\otimes A}$ and
${A\otimes 1}$ are subsets of the closed linear spans of
${(A\otimes 1)\DD A}$ and ${(1\otimes A)\DD A}$, respectively
and therefore each of these closed linear spans is equal to
${A\otimes A}$, thereby ensuring $\AA$ is a compact quantum
group. Co-amenability is then immediate. We shall show only
that $L=A$; the proof that $R=A$ is similar. Arguing by
contradiction, suppose that $L\ne A$, so that there exists a
non-zero element $\t\in A^*$ that vanishes on~$L$. Then
$\t(a*hb)=0$, for all $a,b\in A$. Thus, ${(\t\otimes h
b)\DD(a)=0}$; that is, $h(b((\t\otimes \id)\DD(a)))=0$. By
faithfulness of $h$ we deduce that $(\t\otimes \id)\DD(a)=0$.
Applying $\e$ now we get $0= \e((\t\otimes
\id)\DD(a))=\t((\id\otimes \e)\DD(a))=\t(a)$. Hence, $\t=0$, a
contradiction. Therefore, $L=A$, as required.~\qed

The question arises as to whether one can drop the
faithfulness requirement on the Haar integral $h$ in the
preceding theorem. The answer is no. To see this let
$A=C({\mathbf D})$, the C*-algebra of continuous
complex-valued functions on the closed unit disc ${\mathbf
D}$. A co-multiplication $\DD$ on $A$ is given by setting
$\DD(f)(s,t)=f(st)$, for all $s,t\in \bC$. The linear
functional $\d_0$ on $A$ defined by evaluation at the origin,
$\d_0(f)=f(0)$, is a Haar integral for $\AA$ and the
functional $\d_1$ is a bounded co-unit. But $\AA$ is not a
compact quantum group, by \cite[Proposition~2.2]{MT}.

\section{Appendix}

For the convenience of the reader we gather here some basic
facts about compact quantum groups
( see \cite{KT,MT,Wo2} for more information).

 A compact quantum group $(A, \Delta )$ consists of
 a unital $C*$-algebra $A$ and a unital $*$-homomorphism
 $\DD : A \rightarrow A \otimes A$ (called the co-multiplication)
 satisfying
\[ (\id \otimes \DD) \DD = (\DD \otimes \id) \DD \]
 and such that the linear spans of $(1 \otimes A) \DD A$ and
 $(A \otimes 1) \DD A$ are each dense in $A\otimes A$. A
{\em morphism} from $\AA$ to a compact quantum group $(B,\DD')$ is
a unital $*$-homomorphism ${\pi\colon A\to B}$ satisfying
$\DD' \pi={(\pi\otimes \pi)\DD}$.

There exists a unique state $h$ on $A$ called the Haar integral
 of $(A, \Delta )$ which satisfies

  \[ (h \otimes \id) \DD = (\id \otimes h) \DD =  h(\cdot) 1. \]

By a Hopf $*$-subalgebra $\A$ of a compact quantum group
$(A, \DD )$ we mean a Hopf $*$-algebra such that $\A$ is a
$*$-subalgebra of~$A$ with co-multiplication given by 
restricting the co-multiplication $\DD$ from $A$ to $\A$.
The co-unit $\varepsilon \colon  \A \rightarrow \bC $ and the
co-inverse  $\kappa \colon  \A \rightarrow \A$ of $\A$ are
linear maps satisfying

\[ (\e \otimes \id)\DD = (\id \otimes \e)\DD = \id, \]
\vspace{-2ex}
\[ m (\kappa \otimes \id) \DD = m(\id \otimes \kappa)\DD =  \e(\cdot) 1, \]

\noindent where $m \colon \A \otimes \A \rightarrow \A$
denotes the multiplication map. The co-unit $\e$ is known to be a
$*$-homomorphism.

\vspace{1ex}
Any compact quantum group $(A, \DD )$ has a canonical dense Hopf
$*$-subalgebra $\A$ consisting of the linear span
of the matrix entries of all finite dimensional co-representations
of $(A, \DD )$.
By abuse of language  $\e$ and $\kappa$  are also refered
to as the co-unit and the co-inverse of $(A, \DD )$. We call $\A$
the associated Hopf $*$-algebra of $(A, \DD )$.

\vspace{2ex}
The associated Hopf $*$-algebra of a compact quantum group
has the following uniqueness property (which is stated without proof
in \cite{KT}).

\begin{thm}
The associated Hopf $*$-algebra $\A$ of a compact quantum
group $(A, \DD )$ is the unique dense Hopf $*$-subalgebra of $(A, \DD )$.
\end{thm}

\demo  Let $\B$ be another dense Hopf $*$-subalgebra of $(A, \DD )$.
We must show that $\A =\B$.
First we show that $\B$ is the linear span of the matrix entries of
those finite-dimensional
co-representations which have matrix entries belonging to $\B$.
This will immediately imply that $\B\subset\A$.
Thus let $x\in\B$. Then we may write $\DD (x)=\sum_i x_i\otimes y_i$,
for finitely many  $x_i, y_i\in\B$ with $\{y_i\}$
linearly independent. Pick linear functionals $\{\xi_i\}$ on $\B$ such
that
$\xi_i (y_j )=\delta_{ij}$, for all $i,j$.
Then
$$\DD (x_i )=(\id\otimes\id\otimes\xi_i )\sum_j\DD (x_j )\otimes y_j
=(\id\otimes\id\otimes\xi_i )(\DD\otimes\id )\DD (x)$$
$$=(\id\otimes\id\otimes\xi_i )(\id\otimes\DD )\DD (x)
=\sum_j x_j\otimes (\id\otimes\xi_i )\DD (y_j ),$$
for all $i$. Thus, if we let $\{e_i\}$ denote a linear basis for the
vector subspace of $\B$ spanned by
$\{x_i\}$, there exist finitely many elements $z_i, w_{kl}$ in $\B$
such that
$$\DD (x)=\sum_i e_i\otimes z_i\ \ \ \ \ {\mathrm{and}}\ \ \ \ \ \DD
(e_j )=\sum_k e_k\otimes w_{kj}\ , $$
for all $j$. Now
$$\sum_{k,l}e_l\otimes w_{lk}\otimes w_{kj}=\sum_k\DD (e_k )\otimes
w_{kj}=(\DD\otimes\id )\DD (e_j )$$
$$=(\id\otimes\DD )\DD (e_j )=\sum_l e_l\otimes\DD (w_{lj} ),$$
so by linear independence of $\{e_i\}$, we get
$$\DD (w_{lj})=\sum_k w_{lk}\otimes w_{kj}\ ,$$
for all $j,l$. It follows that $w=(w_{ij})$ is a
finite-dimensional
co-representation of $(A, \DD )$ with matrix entries belonging to
$\B$. Furthermore, the element $x$
is a linear combination of the matrix entries of $w$ because
$$x=(\id\otimes\varepsilon )\DD (x)=\sum_i \varepsilon (z_i )e_i
=\sum _i \varepsilon (z_i )(\varepsilon\otimes\id )\DD (e_i
)=\sum_{ij}\varepsilon (z_i e_j )w_{ji}\ ,$$
where $\varepsilon$ is the co-unit of $\B$. This proves that
$\B\subset\A$.

To prove the converse inclusion,
first observe that $\B$ is the linear span of the matrix
entries of those finite-dimensional irreducible unitary
co-representa- \\ tions of $(A, \DD )$ with matrix entries
belonging to $\B$. To see this, consider the co-representation
$w$ constructed above, and define elements
$v_{ij}=w_{ij}+(\delta_{ij}-\varepsilon (w_{ij}))I\in\B$,
for all $i,j$, where $\varepsilon$ now denotes the co-unit of $\A$.
It is easily checked that $v=(v_{ij})$ is a
co-representation of $(A, \DD )$.
Since $\varepsilon (v_{ij})=\delta_{ij}$, for all $i,j$,
the co-representation $v$ is invertible
with inverse $v^{-1}=(\kappa (v_{ij}))$, where $\kappa$ is the  co-inverse
of $\A$.
Now it is known   \cite{ KT, Wo1} that any invertible
co-representation is equivalent to a direct sum of
irreducible unitary ones. Since the invertible co-representation $v$
has matrix entries in $\B$, its irreducible components are easily
seen to also have matrix entries belonging to $\B$. It then follows that
$\B$ is a linear span of the required sort.

To conclude that $\A \subset \B$, we  now show that every
finite-dimensional irreducible unitary co-representation
of $(A,\DD)$ is equivalent to one with matrix entries belonging to
$\B$. Assume, for contradiction, that $v =(v_{ij})$ is
a finite-dimensional irreducicble unitary co-representation
not equivalent to any finite-dimensional irreducicble unitary
co-representation $u=(u_{ij})$ with matrix entries $u_{ij}$
belonging to $\B$. From \cite[Theorem 7.4]{MT}, we get that
$h(u_{ij}v_{kl} )=0$, for all $i,j,k,l$. Since
$\B$ is linearly spanned by elements of the
type $u_{ij}$, as observed above, and $\B$ is dense in $A$, this
implies that
$h(a v_{kl} )$, for all $k,l$ and $a\in A$.
In particular, we get $h(v_{kl}^* v_{kl})=0$, and therefore $v_{kl}=0$,
for all $k,l$, since $h$ is faithful on $\A$. This is impossible
as $v$ is unitary.~\qed

Note that the first part of the proof shows that $\A$ is
maximal among all Hopf $*$-subalgebras of $(A, \DD )$.

\bigskip
{\parindent=0pt Addresses of the authors:

\smallskip Erik B\'edos, Institute of Mathematics, University of 
Oslo, \\ 
P.B. 1053 Blindern, 0316 Oslo, Norway. E-mail: bedos@math.uio.no. \\

Gerard J. Murphy, Department of Mathematics, National University of
Ireland,
Cork, Ireland. E-mail: gjm@ucc.ie.\\

Lars Tuset, Faculty of Engineering, Oslo University College, \\ 
Cort Adelers Gate~30, 0254 Oslo, Norway. E-mail:Lars.Tuset@iu.hio.no.}

\end{document}